\documentclass[12pt]{article}

\usepackage{amsfonts,amsmath,amssymb,comment}
\usepackage[english]{babel}
\usepackage[usenames]{color}
\usepackage{setspace}

\RequirePackage{natbib}
\RequirePackage[colorlinks,citecolor=blue,linkcolor=blue,urlcolor=blue,pagebackref]{hyperref}
\RequirePackage{hypernat}
\hypersetup{breaklinks=true}

\newcommand{\be}{\begin{equation}}
	\newcommand{\ee}{\end{equation}}
\newcommand{\bea}{\begin{eqnarray}}
	\newcommand{\eea}{\end{eqnarray}}
\newcommand{\beas}{\begin{eqnarray*}}
	\newcommand{\eeas}{\end{eqnarray*}}

\newtheorem{theorem}{Theorem}[section]
\newtheorem{definition}[theorem]{Definition}
\newtheorem{proposition}[theorem]{Proposition}

\newtheorem{lemma}[theorem]{Lemma}
\newtheorem{remark}[theorem]{Remark}
\newtheorem{example}[theorem]{Example}
\newtheorem{examples}[theorem]{Examples}
\newtheorem{foo}[theorem]{Remarks}
\newenvironment{Example}{\begin{example}\rm}{\end{example}}

\newenvironment{proofwithoutproofatstart}{\noindent}
{\unskip\nobreak\hfill$\Box$\par}

\newcommand{\E}{\mathbb{E}}

\DeclareMathOperator{\argmin}{arg\,min}

\DeclareMathOperator{\Imp}{\text{Im}(\phi)}

\newcommand{\real}{\mathbb{R}}

\def\<{\langle}
\def\>{\rangle}

\def\Re{\mathop{\operator@font Re}}
\def\Im{\mathop{\operator@font Im}}
\def\Erl{\mathop{\operator@font Erl}}
\def\Det{\mathop{\operator@font Det}}
\def\eq{\mathop{\operator@font eq}}

\def\a{{\alpha}}

\def\l{{\lambda}}

\def\<{\langle}
\def\>{\rangle}

\begin{document}
	
	\title{\vskip -2.2cm A Rank-Dependent 
		Theory \\ for Decision under Risk and Ambiguity\thanks{We are very grateful to Patrick Cheridito, Louis Eeckhoudt,
			Itzhak Gilboa, Glenn Harrison,
			Edi Karni,
			Mark Machina, John Quiggin, Frank Riedel, Harris Schlesinger,
			Peter Wakker,
			and seminar and conference participants at EURANDOM, Leibniz University Hannover, Princeton University, Tilburg University,
			the Tinbergen Institute, the University of Amsterdam, the European Meeting of Statisticians, 
			the CEAR RDU Workshop 
			in Atlanta,
			and the Bielefeld-Hannover Workshop on Risk Measures and Uncertainty in Insurance for their comments and suggestions.
			An earlier version of this paper was circulated under the title
			``A dual theory for decision under risk and ambiguity''.
			This research was funded in part by the Netherlands Organization for Scientific Research
			under grants NWO-VIDI and NWO-VICI 2020--2027 (Laeven) and grant NWO-VENI (Stadje).}}
	
	\author{Roger J.~A.~Laeven\\
		{\footnotesize Department of Quantitative Economics}\\
		{\footnotesize University of Amsterdam,}\\
		{\footnotesize EURANDOM and CentER}\\
		{\footnotesize {\tt R.J.A.Laeven@uva.nl}}\\
		\and Mitja Stadje\\
		{\footnotesize Faculty of Mathematics and Economics}\\
		{\footnotesize Ulm University}\\
		{\footnotesize {\tt Mitja.Stadje@uni-ulm.de}}}
	
	\date{
		This Version: \today}\maketitle
	\vspace{-0.9cm}
	\begin{abstract}
		This paper axiomatizes,
		in a two-stage setup,
		a new theory for decision under risk and ambiguity.
		The axiomatized preference relation $\succeq$ on the space $\tilde{V}$ of random variables induces
		an ambiguity index $c$ on the space $\Delta$ of probabilities, 
		a probability weighting function $\psi$,
		generating the measure $\nu_{\psi}$ by transforming an objective probability measure,
		and a utility function $\phi$,
		such that, for all $\tilde{v},\tilde{u}\in\tilde{V}$,
		\begin{align*} 
			\tilde{v}\succeq\tilde{u}
			\Leftrightarrow &\min_{Q \in \Delta}
			\left\{\E_Q\left[\int\phi\left(\tilde{v}^{\centerdot}\right)\,\mathrm{d}\nu_{\psi}\right]+c(Q)\right\}\geq
			\min_{Q \in \Delta}
			\left\{\E_Q\left[\int\phi\left(\tilde{u}^{\centerdot}\right)\,\mathrm{d}\nu_{\psi}\right]+c(Q)\right\}.
		\end{align*}
		Our theory extends the rank-dependent utility model of \cite{Quiggin82}
		for decision under risk
		to risk \textit{and} ambiguity,
		reduces to the variational preferences model 
		when $\psi$ is the identity,
		and is \textit{dual} to variational preferences when $\phi$ is affine
		in the same way as the theory of \cite{Yaari87} is dual to expected utility. 
		As a special case, we obtain a preference axiomatization of a decision theory that is
		a rank-dependent generalization of
		the popular maxmin expected utility theory. 
		We 
		characterize ambiguity aversion in our theory. 
	\end{abstract}
    
	
	\newpage \onehalfspacing
	
	\setcounter{equation}{0}
	
	\section{Introduction}
	
	The distinction between risk (probabilities given) and ambiguity (probabilities unknown),
	after \cite{Keynes21} and \cite{Knight21},
	has become a central aspect in decision-making under uncertainty.
	Already since \cite{Ellsberg61} the importance of this distinction had been apparent:
	whereas in the classical subjective expected utility (SEU) model of \cite{Savage54}
	the distinction between risk and ambiguity was nullified
	through the assignment of subjective probabilities (\citealp{Ramsey31}, \citealp{deFinetti31}),
	the \cite{Ellsberg61} paradox showed experimentally that decisions under ambiguity
	could not be reconciled with any such assignment of subjective probabilities.
	It took, however, until the 1980s before decision models were developed that
	could account for ambiguity without the assignment of subjective probabilities.
	
	Among the most popular---by now canonical---models for decision under risk and ambiguity today are
	maxmin expected utility (MEU, \citealp{GilboaSchmeidler89}), also called multiple priors,
	and Choquet expected utility (CEU, \citealp{Schmeidler86, Schmeidler89}).
	The former model is a decision-theoretic foundation of the classical decision rule of \cite{Wald50} in (robust) statistics;
	see also \cite{Huber81}.
	Somewhat more recently, \cite{MMR06} axiomatized the broad and appealing class of variational preferences (VP),
	which includes MEU and the multiplier preferences of \cite{HansenSargent00, HS01} as special cases.
	Multiplier preferences have been widely used in macroeconomics,
	to achieve ``robustness'' in settings featuring model uncertainty.
	
	In the \cite{AnscombeAumann63} setting, all the aforementioned decision models reduce
	to the classical \cite{NeumannMorgenstern44} expected utility (EU) model under risk---a property that is undesirable from a descriptive perspective:
	it means, for example, that the \cite{Allais53} paradox, hence the common consequence and common ratio effects,
	are still present under risk; see e.g., \cite{Machina87}.\footnote{Furthermore, \cite{Machina09} shows that decision problems in the style of \cite{Ellsberg61} lead to similar paradoxes for CEU as for SEU,
		arising from event-separability properties that CEU retains in part from SEU.} 
	
	In this paper, we introduce and axiomatize,
	in a two-stage setup similar to 
	the \cite{AnscombeAumann63} setting,
	a new theory for decision under risk and ambiguity.
	As we will explicate, our theory extends the rank-dependent utility (RDU) model of \cite{Quiggin82}
	for decision under risk to risk \textit{and} ambiguity;
	reduces to the VP model of \cite{MMR06}
	when linear in objective probabilities;
	and is \textit{dual} to VP when affine in wealth
	in the same way as the theory of \cite{Yaari87}
	is dual to EU for decision under risk.
	Thus, the theory developed in this paper may be viewed as an extension of the dual theory (DT) of \cite{Yaari87}
	and the RDU model of \cite{Quiggin82}
	for decision under risk to settings involving risk and ambiguity,
	just like the theory of \cite{MMR06} is a significant extension to risk and ambiguity of the EU model for risk.
	As a special case, we obtain a preference axiomatization of a decision
	theory that is a rank-dependent generalization of
	the popular MEU model of \cite{GilboaSchmeidler89},
	and is dual to it when affine in wealth.
	See Table~\ref{tab:dual}.\footnote{Our approach can also be used to obtain decision theories that are dual to, and a rank-dependent generalization of,
		CEU of \cite{Schmeidler89}. 
		We do not pursue this in this paper.
		Similar comments apply to the increasingly popular $\alpha$-maxmin expected utility model and the special case of the Hurwicz expected utility model (\citealp{GMM04}, \citealp{GulPesendorfer15}).}
	\begin{table}[h]\label{tab:dual}
		\caption{Primal, Dual and Rank-Dependent Decision Theories}
		\centering
		\vskip 0.1cm 
		\begin{tabular}{ c | c | c | c }
			& \textit{Primal} & \textit{Dual} & \textit{Rank-Dependent}\\
			\hline\hline
			Risk & EU (vNM, 1944) & DT (Yaari, 1987) & RDU (Quiggin, 1982)\\
			\hline
			Risk and ambiguity & MEU (GS, 1989) & \textbf{\textit{This paper}} & \textbf{\textit{This paper}}\\
			\hline
			Risk and ambiguity & VP (MMR, 2006) & \textbf{\textit{This paper}} & \textbf{\textit{This paper}}\\
			\hline
		\end{tabular}
	\end{table}
	
	The development of the DT 
	of \cite{Yaari87}
	was methodologically motivated by the fact that,
	under EU, the decision-maker's (DM's) attitude towards wealth,
	as represented by the utility function,
	completely dictates the attitude towards risk.
	However, attitude towards wealth and attitude towards risk should arguably be treated separately:
	they are ``horses of different colors'' \citep{Yaari87}.
	This is achieved within the DT and RDU models.
	From an empirical perspective, the DT and RDU models naturally rationalize
	various behavioral patterns that are inconsistent with EU. 
	The RDU model synthesizes and encompasses both EU and DT,
	and serves as the main building block in prospect theory of \cite{TverskyKahneman92} discussed in much detail in \cite{Wakker10}.\footnote{According to \cite{HarrisonSwarthout16}, RDU even arises as the most important non-EU model for decision under risk
		from a descriptive perspective.}$^{,}$\footnote{
		Contrary to the linearity in probabilities that occurs under EU, to which MEU and VP reduce under risk,
		the DT model of \cite{Yaari87} is affine in wealth.
		\cite{Yaari87} suggests the behavior of a profit maximizing firm as a prime example in which affineness in wealth seems particularly suitable.
		Other theories that stipulate affineness in wealth
		are provided by convex measures of risk
		(\citealp{FollmerSchied16}, Chapter 4) 
		encompassing many classical insurance premium principles, 
		and by robust expectations (see e.g., \citealp{Riedel09}, and the references therein).
		Despite the popularity of these theories,
		neither affineness in wealth nor linearity in probabilities as in EU is considered fully empirically viable
		for individual decision-making.
		Instead, we provide a general decision theory for risk and ambiguity in which preferences under risk are represented by
		the more general measure on the wealth-probability plane given by RDU.
	}
	
	Similarly, our results are both theoretically (i.e., methodologically) relevant
	and of potential empirical interest.
	At the methodological level, our theory disentangles attitude towards wealth from
	attitude towards risk \textit{and} attitude towards ambiguity.
	We characterize (comparative) ambiguity aversion in our decision model to
	corroborate this separation.
	From an empirical perspective,
	an important and distinctive feature of our theory
	is that it accounts for ambiguity---hence is not subject to violations of subjective probabilities such as the Ellsberg paradox---and yet does not collapse to EU under risk
	as would be the case if the \cite{AnscombeAumann63} setting would apply---hence is not subject to the objective phenomena of the Allais paradox and related effects. 
	We offer three additional motivations for our new decision theory in Section~\ref{sec:motivation}.
	
	
	The numerical representation of the decision theory we axiomatize entails that the DM considers,
	for each random variable to be evaluated in the face of risk and ambiguity,
	a collection of potential probabilistic models rather than a single probabilistic model.
	In recent years, we have seen increasing interest in optimization, macroeconomics, finance and other fields
	to account for the possibility that an adopted probabilistic model
	is an approximation to the true probabilistic model and may be misspecified.
	Models that explicitly recognize potential misspecification provide a ``robust'' approach.
	Within the MEU model, the DM assigns the same plausibility to each probabilistic model
	in a set of probabilistic models under consideration.
	The multiplicity of the set of probabilistic models then reflects the degree of ambiguity.
	The VP model significantly generalizes the MEU model
	by allowing to attach a plausibility (or ambiguity) index to each probabilistic model.
	Such an ambiguity index also appears in the numerical representation of our decision model.
	
	More specifically, our numerical representation $U$ of the preference relation $\succeq$
	on the space $\tilde{V}$ of random variables takes the form
	\begin{equation}\label{eq:numrep}
		U(\tilde{v})=\min_{Q \in \Delta}\left\{\E_Q\left[\int\phi\left(\tilde{v}^{\centerdot}\right)\,\mathrm{d}\nu_{\psi}\right]+c(Q)\right\},\qquad \tilde{v}\in\tilde{V},
	\end{equation}
	with $\Delta$ a set of probabilities on the states of the world,
	$c:\Delta\rightarrow[0,\infty]$ the ambiguity index,
	$\psi:[0,1]\rightarrow[0,1]$ a probability weighting function,
	$\nu_{\psi}$ a measure obtained by transforming an objective probability measure according to $\psi$,
	and $\phi:\mathbb{R}\rightarrow\mathbb{R}$ a utility function.
	Special cases of interest occur when the ambiguity index is the well-known relative entropy or Kullback-Leibler divergence (\citealp{HansenSargent00, HS01} and \citealp{Strzalecki11a}),
	or, more generally, an $f$-divergence measure
	(\citealp{Csiszar75}, \citealp{Ben-Tal85}, and \citealp{LaevenStadje13}),
	or simply an indicator function that takes the value zero on a subset of $\Delta$ and $\infty$ otherwise.
	It is directly apparent from \eqref{eq:numrep} that
	in the absence of uncertainty about the state of the world (i.e., in the case of risk)
	our decision model reduces to RDU of \cite{Quiggin82}.
	The familiar utility and probability weighting functions determine the attitude towards wealth and risk.
	In our general model, we characterize ambiguity  aversion in terms of the ambiguity index.
In essence, our axiomatization is based on a modification of two axioms
stipulated by \cite{GilboaSchmeidler89} and \cite{MMR06}:
the uncertainty aversion axiom and the (weak) certainty independence axiom. 
First, we postulate a form of ambiguity aversion (Axiom A6 below)
with respect to ``subjective mixtures of random variables''
rather than with respect to ``probabilistic mixtures of horse lotteries'' as in the uncertainty aversion axiom of \cite{GilboaSchmeidler89} and \cite{MMR06}.
Subjective mixtures of \textit{horse lotteries} are due to \cite{Nakamura90} and \cite{Gul92}
(see also \citealp{GhirardatoMarinacci01}),
within a (very) different setting and for a different purpose.
The subjective mixtures are employed here for \textit{random variables}, and are extended to subjective additions of random variables.\footnote{Subjective mixtures
	of horse lotteries entail that
	the utility profile of the subjective mixture of two horse lotteries
	equals the convex combination of the utility profiles of the horse lotteries themselves.
	Similarly, the utility profile of a subjective mixture 
	of two random variables
	equals the convex combination 
	of the utility profiles of the random variables themselves.}

Consider two random variables with an unknown probability distribution between which the DM is indifferent.
Then, our new Axiom A6 stipulates that the DM prefers 
a subjective mixture of the two random variables
to 
either one in full.
This constitutes a preference for \textit{diversification},
induced by subjective mixtures of random variables with an unknown probability distribution.\footnote{Under affineness in wealth, this reduces to a preference for convex combinations of the two random variables.}
In the primal theories of \cite{GilboaSchmeidler89} and \cite{MMR06},
uncertainty aversion instead takes the form of a preference for \textit{randomization}.
Translated to our setting of preferences over random variables,
randomization stipulates that the DM prefers receiving two 
random variables, between which she is indifferent,
with probabilities $p$ and $1-p$, $0<p<1$, to obtaining one of them with certainty.
Randomization as used in the primal theories provides a hedge against ambiguity
by trading off ambiguity for chance; diversification as used in our theory provides a hedge against ambiguity
by securing (utility units of) wealth.

Second, we replace the (weak) certainty independence axiom of \cite{GilboaSchmeidler89} and \cite{MMR06}
by a comonotonic type of independence axiom (Axiom A7 below),
which pertains to subjective mixtures and subjective additions of random variables instead of probabilistic mixtures of 
lotteries.\footnote{As is well-known, the independence axiom and its various alternatives
	are key to obtaining, and empirically verifying, preference representations.}
Our approach is inspired by the ``dual independence axiom'' 
of \cite{Yaari87}.
However, in our general setting that allows for a set of probabilistic models, 
the implications of comonotonicity, and its interplay with ambiguity, must be reconsidered:
whereas preferences over random variables may well be invariant to the (subjective) addition
of comonotonic random variables when probability distributions are given
(i.e., under risk)
as 
stipulated by \cite{Yaari87} and \cite{Quiggin82},
this implication may no longer be appropriate under ambiguity
because such addition may impact the ``level'' of ambiguity (see Example~\ref{ex:2} below).

Therefore, we postulate the following two versions of the dual independence axiom
to extend RDU to a setting featuring risk and ambiguity:
($i$) preferences over random variables are invariant to the subjective addition of a comonotonic random variable
with an objectively given probability distribution (Axiom A7);
and ($ii$) preferences over random variables are invariant to subjective mixtures of the random variable and a
comonotonic random variable with an objectively given probability distribution (Axiom A7$^{0}$).
The former yields a decision theory that is a rank-dependent generalization of VP,
the latter yields a decision theory that is a rank-dependent generalization of MEU.
The mathematical details in the proofs of our characterization results are delicate.

\subsection{Related Literature}
In interesting and important work, \cite{DeanOrtoleva17} axiomatize a decision theory that, like our theory,
simultaneously allows for both violations of EU under risk \`a la Allais
as well as violations of SEU \`a la Ellsberg.
Their novel preference representation takes the form of a maxmin multiple priors-multiple weighting functional
that distorts objective probabilities by finding the worst from a class of probability weighting functions.
We briefly highlight the main differences between their work and ours.
Our representation allows for a non-trivial ambiguity index,
including e.g., Hansen-Sargent type robustness via Kullback-Leibler divergences, or general $f$-divergences.
Furthermore, whereas the elements in the class of probability weighting functions in \cite{DeanOrtoleva17} are all increasing and convex,
our theory allows for a single, general probability weighting function,
including star-shaped or (inverse) S-shaped functions often found in experiments (see e.g., \citealp{TverskyKahneman92}, \citealp{Prelec98} and \citealp{GonzalezWu99}) and non-decreasing functions consistent with Value-at-Risk and Expected Shortfall measures of risk (see Section~\ref{sec:robustrm}).
From a methodological perspective, we introduce a setting where, in the spirit of \cite{Yaari87}, random variables, instead of sets of lotteries and acts as in \cite{DeanOrtoleva17},
are the basic objects.
This not only fundamentally impacts the meaning of some of the axioms compared to the setting in \cite{DeanOrtoleva17},
but may also facilitate their interpretation.
In particular, we do not stipulate preference relationships over certain \textit{sets} of lotteries and acts
but consider simple preferences over (single)
random variables.
Finally, our setting enables us to identify and exploit (in our technical proofs) dual relations between the primal theories given by VP and MEU and our theory (see also Table~\ref{tab:dual}).\footnote{More specifically, our proofs employ the (`dual') space of conditional reflected quantile functions and  convex duality results for niveloids.}
These dual relations constitute one of the main pillars underlying the results in this paper.

%
%
%
%
%
%

In an early version of this paper,\footnote{Available from \url{https://pure.tue.nl/ws/files/32872072/024\_report.pdf}.} we extended Yaari’s DT to risk \textit{and} ambiguity.
Then, subjective mixtures of random variables are not required; it is sufficient to consider convex combinations of the random variables themselves.
The (substantially) generalized setting in the current version of the paper allows for utility functions that are non-affine in wealth, extending RDU of \cite{Quiggin82} to risk \textit{and} ambiguity.  
In concurrent work, \cite{Wang22} adopts a two-stage setup that maintains the MEU (or VP) framework in the first stage and replaces EU under risk 
by an RDU representation in the second stage on a finite objective probability space.
Whereas the numerical representation obtained in \cite{Wang22} is similar to that in the present version of this paper, the two approaches are (very) different, both methodologically and technically.
In \cite{Wang22}, the axioms of \cite{GilboaSchmeidler89} (or \citealp{MMR06}) over acts are maintained in the first stage. 
The RDU model is then induced directly in the second stage, by introducing an axiom referred to as p-trade-off consistency, which postulates that trade-offs between objective probabilities of certain ordered outcomes can be done independently of the other outcomes. 
Furthermore, the usual mixtures of (objective) lotteries are replaced by subjective mixtures of \emph{probabilities}. 
In our approach, instead, random variables are the basic objects, leveraging duality in the sense of \cite{Yaari87}. 
Our approach introduces axioms that are different, in meaning and interpretation, from those in \cite{GilboaSchmeidler89} (or \citealp{MMR06}), involving subjective mixtures of random variables and an associated diversification preference and dual independence condition, operating on a rich probability space.  
Thus, the two approaches rely on different behavioral axioms with their own basic objects, mixtures and spaces, and corresponding proof techniques. 

\subsection{Outline}
This paper is organized as follows.
In Section~\ref{sec:setup}, we introduce our setting and notation.
In Section~\ref{sec:rep}, we review some preliminaries, introduce subjective mixtures of random variables and our new axioms, state our main representation results, and discuss their interpretation.
Section~\ref{sec:ambrisk} characterizes ambiguity  aversion in our theory.
In Section~\ref{sec:motivation}, we further motivate our new decision theory from three different perspectives.
All proofs are relegated to the Appendix.

\setcounter{equation}{0}

\section{Setup and Notation}\label{sec:setup}

We adopt a two-stage setup as in the \citet{AnscombeAumann63} approach,
with a state space that is a product space
admitting a two-stage decomposition.
Different from the \citet{AnscombeAumann63} model, however,
we do not assume nor induce EU for risk.
Our theory defines a preference relation over random variables
just like \citet{Yaari87} for risk by which our notation is inspired.
We now formalize our two-stage setup in detail.

We consider a possibly infinite set $W$ of states of the world with $\sigma$-algebra $\Sigma'$ of subsets of $W$ that are events. 
We call a function $F:W\times \mathbb{R}\rightarrow [0,1]$
a conditional cumulative distribution function (CDF)
if, for all $w\in W$, $F(w,\cdot)$ is a CDF and, for every $t\in\mathbb{R}$, $F(\cdot, t)$ is
$\Sigma'$-measurable. 
%
Furthermore, consider, in neo-Bayesian nomenclature, an \textit{act} or \textit{horse (race) lottery} 
$f:W\rightarrow \mathcal{X}$, where $\mathcal{X}$ is an affine space of consequences.
Assume in particular that $\mathcal{X}$ is 
given by the space of \textit{objective 
lotteries}, i.e., the space of CDFs with bounded support.
Then, our preference domain of random variables 
can be seen to correspond and be equivalent (up to a richness condition) to that of the standard Anscombe-Aumann setting as in \cite{GilboaSchmeidler89} and \cite{MMR06}, where the primitives are acts with consequences that are objective lotteries. 
More specifically, as we will see below, we can formally identify random variables and associated conditional CDFs with acts. 
    
In principle, we could---as in much of the modern decision-theoretic literature---formulate our axioms below in terms of acts.
Instead, we define our axioms in terms of random variables, for clarity of exposition and ease of interpretation.
The essential difference compared to \cite{GilboaSchmeidler89} and \cite{MMR06} occurs in our Axioms A6, A7 and A8 below.
Here, among other aspects, the ``$+$'' (and ``$\oplus$'') operation is defined \emph{outcome-wise} instead of \emph{probability-wise}. 
More specifically,  
in \cite{GilboaSchmeidler89} and \cite{MMR06}, the numerical representation of preferences is linear in $\mathcal{X}$ and takes the EU form under risk. 
Formally, in their setting, for two consequences $x_1,x_2\in \mathcal{X}$, the mixing operation 
$	\alpha x_1+(1-\alpha)x_2$ 
represents the compound lottery that yields consequence $x_1$ with probability $\alpha$ and consequence $x_2$ with probability $(1-\alpha)$. 
By contrast, in our axioms, 
the ``$+$'' operation will refer to mixing the outcomes (i.e., payoffs) of the consequences $x_1,x_2$, meaning that, for $\alpha x_1+(1-\alpha)x_2$, the DM receives $100 \cdot \alpha\%$ of the payoff of the consequence $x_1$ and $100 \cdot (1-\alpha)\%$ of the payoff of the consequence $x_2$. 
Furthermore, in our setting, mixing of payoffs will only be considered for comonotonic consequences, which move in tandem.
This is most clearly expressed in terms of random variables.\footnote{ 
In the setting of decision under risk, this 
corresponds to the dual theory of \citet{Yaari87}. 
}
    
Consider next a non-atomic probability space $(S,\Sigma,P)$, supporting a random variable $U$ that is uniformly distributed on the unit interval under $P$.
The product space of interest is given by $W\times S$.
Assume that for every $A\in\Sigma$, the mapping $w\mapsto P[A]$ is $\Sigma'$-measurable.
Let $\tilde{V}$ be the space of all bounded random variables
$\tilde{v}$ defined on the space $(W\times S,\Sigma'\otimes\Sigma)$, i.e.,
$\tilde{v}$ is a mapping from $W\times S$ to a bounded subset of $\mathbb{R}$.
Similar to \citet{Yaari87}, realizations of the random variables in $\tilde{V}$ will be
viewed as payments denominated in monetary units.
For a random variable $\tilde{v}\in\tilde{V}$ and a fixed $w\in W$, the random variable $\tilde{v}^w:S\rightarrow \mathbb{R}$ is the
outcome $s$ contingent payment that the DM receives if she lives in state of the world $w$.
This makes $\tilde{v}^w$ also interpretable as a \textit{roulette lottery}, in neo-Bayesian nomenclature.
Henceforth, $\tilde{v}^w$ and its associated roulette lottery are often identified.\footnote{We have assumed here that in each state of the world $w$ the possible outcomes are the same.
This can simply be achieved by adding to each state of the world additional outcomes with associated probability zero.
}
We denote by $\tilde{V}_0$ the subspace of all random variables in $\tilde{V}$ that take only finitely many values. 
For fixed $w\in W$, we define 
the conditional CDF
$F_{\tilde{v}}(w,t)$ of the $\Sigma$-measurable random variable
$\tilde{v}^{w}$, given by $s\mapsto\tilde{v}^{w}(s)$, by
$F_{\tilde{v}}(w,t)=P[\tilde{v}^{w}\leq t]$.
(We sometimes omit the dependence on $t$, i.e., we sometimes write $F_{\tilde{v}}(w)$.)
From the assumptions above, it follows that for
every $t\in\real$, $F_{\tilde{v}}(.,t)$ is $\Sigma'$-measurable.

We now identify every act $f:W\rightarrow \mathcal{X}$ with a $\Sigma'\otimes\Sigma$-measurable random variable on the product space $W\times S$, $\tilde{v}: W\times S\rightarrow\mathbb{R}$, in the following way. 
First, every $\tilde{v}\in\tilde{V}$ induces a conditional CDF and hence can be identified with a horse lottery by setting $f(w) := F_{\tilde{v}}(w)$. 
Conversely, for fixed $w$, every horse lottery 
$f$, given by $w\mapsto\mu^w$, for (roulette) lotteries $\mu^w$ defined on $(S,\Sigma,P)$, induces a CDF $F(w,\cdot)$.
Let $q(w)$ be the left-continuous inverse of $F(w)$, i.e.,
	\begin{equation*}
		q(w,\l)=\inf\{t\in\mathbb{R}|F(w,t)\geq \l\},\quad\,\,\,\,\l\in(0,1).
	\end{equation*} 
Then we can define $\tilde{v}^{w}(s)=q(w,U(s))$ and it is easy to see that, for every $w\in W$, $\tilde{v}^{w}$ has the same probability distribution as $f(w)$.
Hence, there is a one-to-one correspondence between equivalence classes of random variables
$\tilde{v}\in\tilde{V}$ with the same conditional distributions and horse lotteries $f$. 
Thus, our preference domain of random variables on the product space $W\times S$ 
corresponds to the standard Anscombe-Aumann domain 
of acts mapping states, $W$, to objective lotteries of outcomes $\mathcal{X}$;
we additionally require richness in the sense of a non-atomic probability space $(S,\Sigma,P)$.
    
	%
%

For some $\tilde{v}\in\tilde{V}$ (e.g., those that represent payoffs
from games such as flipping coins) the DM may
actually know the objective probability distribution.
As in \cite{AnscombeAumann63}, 
\cite{GilboaSchmeidler89} and \cite{MMR06}
the associated objective lotteries will play a special role for our
theory.
In this case, the probability distribution of $\tilde{v}$ does not depend
on $w$, i.e., for all $w_1,w_2\in W$,
$F_{\tilde{v}}(w_1)=F_{\tilde{v}}(w_2)$.
We denote the corresponding space of all
the random variables in $\tilde{V}$ and $\tilde{V}_0$ that
carry no ambiguity by $V$ and $V_0$, respectively.
For random variables $v\in V$ we usually omit the $w$, i.e., we just
write $F_v(t)$ instead of $F_v(w,t)$.
In the space $V$, $v_n$ converges in distribution to $v$ if $F_{v_n}$ converges to $F_v$ for all continuity points
of $F_v$.

Furthermore, let $V'$ be the space defined by\footnote{For instance, for Ellsberg type urns, the space $V'$ corresponds to payoffs depending upon outcomes from urns that contain only balls of one color, but where the specific color is unknown.}
\begin{align*} 
V'=\big\{\tilde{v}\in\tilde{V}&|\tilde{v}\mbox{ is independent of }s\in S, 
\mbox{ i.e., for }s_1,s_2\in S:\,\tilde{v}^{w}(s_1)=\tilde{v}^{w}(s_2)\big\}.
\end{align*}
Clearly, the space $V'$
of all random variables in $\tilde{V}$ that carry no risk
may be identified with the space of bounded measurable functions
on $(W,\Sigma')$.
$V'_0$ is defined as the corresponding subspace of bounded measurable functions that take only finitely many values.

Finally, let $\Delta(W,\Sigma')$ be the space of all finitely additive measures
on $(W,\Sigma')$ with mass one
and denote by $\Delta_{\sigma}(W,\Sigma')$ the
space of all probability measures on $(W,\Sigma')$.

\setcounter{equation}{0}

\section{Representation}\label{sec:rep}

In this section, we provide an axiomatic foundation of a new rank-dependent theory for decision under risk and ambiguity.

\subsection{Preliminaries}
We define a preference relation $\succeq$ on $\tilde{V}_0$.
As usual, $\succ$ stands for strict preference and $\sim$ for indifference.
The preference relation $\succeq$ on $\tilde{V}_0$ induces a preference order,
also denoted by $\succeq$, over random variables $s\mapsto\tilde{v}^{w}(s)$
through those random variables in $\tilde{V}_0$ that are associated with objective lotteries
(i.e., are in $V_0$) by defining
$\tilde{v}^{w}\succeq\tilde{u}^{w}$ if, and only if, $v\succeq u$ with $v,u\in V_{0}$
and $F_{v}(t)=F_{\tilde{v}}(w,t)$, $F_{u}(t)=F_{\tilde{u}}(w,t)$ for all $t\in\mathbb{R}$.
(This, in turn, induces, similarly, a preference relation over monetary payments.)
We assume that $\succeq$ satisfies the following properties:

\vskip 0.14cm
\textit{A{\sc xiom} A1---Weak Order: $\succeq$ is complete and transitive. 
That is:
\begin{itemize}
\item [(a)] $\tilde{v}\succeq\tilde{u}$ or $\tilde{u}\succeq\tilde{v}$ for all $\tilde{v},\tilde{u}\in\tilde{V}_0$.
\item[(b)] If $\tilde{v},\tilde{u},\tilde{r}\in\tilde{V}_0$, $\tilde{v}\succeq\tilde{u}$ and 
$\tilde{u}\succeq\tilde{r}$, then $\tilde{v}\succeq\tilde{r}$.
\end{itemize}
}
\vskip 0.14cm
Whereas \cite{Yaari87} assumes that the preference relation
is complete on the space of all random variables including those taking infinitely many values, 
we only assume in A1(a) that $\succeq$ is complete on
the subspace $\tilde{V}_{0}$ of random variables that take finitely many values.
We will see later that $\succeq$ and our
representation results may 
be uniquely extended to the entire
space $\tilde{V}$.
For the dual theory of~\cite{Yaari87} without the
completeness axiom, see~\cite{Maccheroni04}.
\vskip 0.14cm
\textit{A{\sc xiom} A2---Neutrality: Let $\tilde{v}$ and
$\tilde{u}$ be in $\tilde{V}_0$ and have the same
conditional CDFs, $F_{\tilde{v}}$ and $F_{\tilde{u}}$.
Then, $\tilde{v}\sim\tilde{u}$.}
\vskip 0.14cm
Axiom~A2 states that
$\succeq$ depends only on the conditional distributions.
In particular, $\succeq$ induces a preference relation
on the space of conditional CDFs by defining
$F(\succeq)G$ if, and only if, there exist two random variables
$\tilde{v},\tilde{u}\in\tilde{V}_0$ such that
$\tilde{v}\succeq\tilde{u}$ and $F_{\tilde{v}}(w)=F(w),\ G_{\tilde{u}}(w)=G(w)$ for all $w\in W$.
To simplify notation, we will henceforth use
$\succeq$ both for preferences over random variables and for
preferences over conditional CDFs.
\vskip 0.14cm
\textit{A{\sc xiom} A3---Continuity: For every $v,u\in V_0$ such that
$v\succ u$, and uniformly bounded sequences $v_n$ and $u_n$
converging in distribution to $v$ and $u$, there
exists an $n$ from which onwards $v_n\succ u$ and $v\succ u_n$.
Furthermore, for every $\tilde{v}\in \tilde{V}_0 $, the sets 
$\big\{ m\in\mathbb{R} |m\succ \tilde{v}\big\}
\ \mathrm{and}\ \big\{ m\in\mathbb{R}|\tilde{v}\succ m\big\}\ 
$
are open.}
\vskip 0.14cm
When restricted to $V_0$, our continuity condition~A3 is
equivalent to the one employed in \cite{Yaari87},
where it is a little stronger than that of \cite{MMR06} or \cite{GilboaSchmeidler89}.

\vskip 0.14cm
\textit{A{\sc xiom} A4---Certainty First-Order Stochastic Dominance and Non-Degeneracy: 
For all $v,u\in V_0$:
If $F_v(t)\leq F_u(t)$ for every $t\in\mathbb{R}$, then $v\succeq u$. 
Furthermore, if the inequality is sharp for some $t\in\mathbb{R}$, then $v\succ u$.}
\vskip 0.14cm
\textit{A{\sc xiom} A5---Monotonicity: 
For all $\tilde{v},\tilde{u}\in \tilde{V}_0$:
If $\tilde{v}^w\succeq \tilde{u}^w$ for every $w\in W$, then $\tilde{v}\succeq \tilde{u}$.}
\vskip 0.14cm
We postulate Axioms~A1-A5 for a preference relation $\succeq$ defined on the space of
finite-valued random variables $\tilde{V}_0$.
However, by A2, in view of the one-to-one correspondence explicated in the previous section,
it is straightforward to verify that this preference relation induces a preference relation, also denoted by $\succeq$,
on the space of horse lotteries, satisfying the same axioms.
Consequently, all axioms considered so far, which will be maintained in our setting,
are common; see \cite{Yaari87}, \cite{Schmeidler89}, \cite{GilboaSchmeidler89} and \cite{MMR06}.
To (strictly speaking: a subset of) the collection of axioms above,
\cite{GilboaSchmeidler89} add the following two axioms (involving probabilistic mixtures):
\vskip 0.14cm
\textit{A{\sc xiom} A6MEU---Uncertainty Aversion:
If $\tilde{v},\tilde{u}\in \tilde{V}_0$ and $\alpha\in(0,1)$,
then $F_{\tilde{v}}\sim F_{\tilde{u}}$ implies $\a F_{\tilde{v}}+(1-\a)F_{\tilde{u}}\succeq F_{\tilde{v}}$.}
\vskip 0.14cm
\textit{A{\sc xiom} A7MEU---Certainty Independence: 
If $\tilde{v},\tilde{u}\in \tilde{V}_0$ and $v\in V_0$, then
$F_{\tilde{v}}\succeq F_{\tilde{u}}\,\,\, \Leftrightarrow\,\,\,\a F_{\tilde{v}}+(1-\a)F_{v}\succeq \a F_{\tilde{u}}+(1-\a)F_{v}\quad\mbox{ for all
}\a\in(0,1).
$
}
\vskip 0.14cm
With these axioms at hand, one obtains the maxmin expected utility representation, as follows:\vskip 0.1cm

\noindent{\bf Theorem 4.0.(i) (\citealp{GilboaSchmeidler89})}\\
\textit{A preference relation $\succeq$ satisfies A1-A5 and A6MEU-A7MEU if,
and only if, there exist an increasing and continuous\footnote{\cite{GilboaSchmeidler89} imposed slightly milder continuity and monotonicity conditions than in this paper such that in their setting $\phi$ was non-constant and non-decreasing.} function
$\phi:\mathbb{R}\rightarrow \mathbb{R}$ and a non-empty, closed and
convex set $C\subset \Delta(W,\Sigma')$ such that, for all}
$\tilde{v},\tilde{u}\in\tilde{V}_0$,
\begin{equation*}
\tilde{v}\succeq \tilde{u}\Leftrightarrow
\min_{Q\in C} \E_Q\left[\int \phi(t)F_{\tilde{v}}(.,\,\mathrm{d}t)\right]\geq \min_{Q\in
C}\E_Q\left[\int \phi(t)F_{\tilde{u}}(.,\,\mathrm{d}t)\right].
\end{equation*}
\textit{Furthermore, $\succeq$ has a unique extension to $\tilde{V}$
satisfying the same assumptions (over $\tilde{V}$).}\vskip 0.1cm

Somewhat more recently, \cite{MMR06} obtained a more general representation result, which
includes the maxmin expected utility representation of \cite{GilboaSchmeidler89} as a special case,
but also covers the multiplier preferences employed in robust macroeconomics; see, for
instance, \cite{HansenSargent00, HS01}.
If, in the certainty independence axiom A7MEU, $\alpha$ is close to
zero, then $\a F_{\tilde{v}}+(1-\a)F_v$ carries ``almost no ambiguity''.
Hence, if a DM prefers $\tilde{v}$ to $\tilde{u}$ (as the axiom presumes),
but merely because $\tilde{v}$ carries less ambiguity than $\tilde{u}$,
then she may actually prefer $\a F_{\tilde{u}}+(1-\a)F_{v}$ to $\a F_{\tilde{v}}+(1-\a)F_v$ when $\a$ is small
and ambiguity has almost ceased to be an issue.
Therefore, \cite{MMR06} suggest to replace the certainty independence axiom by the following weaker axiom:
\vskip 0.14 cm
\textit{A{\sc xiom} A7VP---Weak Certainty Independence: 
If $\tilde{v},\tilde{u}\in \tilde{V}_0 ,$ $v,u\in V_0$ and
$\a\in(0,1)$, then}
$\a F_{\tilde{v}}+(1-\a)F_v \succeq \a F_{\tilde{u}}+(1-\a)F_v  \Rightarrow\a F_{\tilde{v}}+(1-\a)F_u\succeq \a F_{\tilde{u}}+(1-\a)F_u.
$
\vskip 0.14 cm
Denote by $m_{\tilde{v}}$ the certainty equivalent of $\tilde{v}$,
that is, $m_{\tilde{v}}\sim \tilde{v}$, $m_{\tilde{v}}\in\mathbb{R}$.
Replacing A7MEU by A7VP (\textit{ceteris paribus}) yields the following theorem:\vskip 0.1cm

\noindent{\bf Theorem 4.0.(ii) (\citealp{MMR06})}\\
\textit{A preference relation $\succeq$ satisfies A1-A5 and A6MEU-A7VP if,
and only if, there exist an increasing and continuous\footnote{\cite{MMR06}, just like \cite{GilboaSchmeidler89}, imposed slightly weaker continuity and monotonicity conditions than in the present paper, 
such that $\phi$ was non-constant and non-decreasing.} function
$\phi:\mathbb{R}\rightarrow \mathbb{R}$ and a
grounded,\footnote{We say that $c$ is \textit{grounded} if $\min_{Q\in
\Delta(W,\Sigma')}c(Q)=0$.} convex and lower-semicontinuous function
$c:\Delta(W,\Sigma')\rightarrow [0,\infty]$ such that, for all
$\tilde{v},\tilde{u}\in\tilde{V}_0$,}
\begin{align*} 
\tilde{v}\succeq
\tilde{u}\,\,\Leftrightarrow\,\, \min_{Q\in \Delta(W,\Sigma')}
\left\{\E_Q\left[\int \phi(t)F_{\tilde{v}}(.,\,\mathrm{d}t)\right]+c(Q)\right\}\qquad\\
\geq\min_{Q\in \Delta(W,\Sigma')}\left\{\E_Q\left[\int
\phi(t)F_{\tilde{u}}(.,\,\mathrm{d}t)\right]+c(Q)\right\}.
\end{align*}
\textit{Furthermore, for each $\phi$ there exists a unique minimal $c_{0}$ given by} \
\begin{equation*}
c_{0}(Q)=\sup_{v'\in V'_0} \left\{m_{\phi(v')}-\E_Q[ \phi(v')]\right\}.
\end{equation*}

\subsection{Subjective Mixtures of Random Variables}\label{sec:submixrv}

In this subsection, we introduce and analyze \textit{subjective mixtures of random variables}.
Let us, following the insightful \cite{GhirardatoMaccheroniMarinacciSiniscalchi03},
first recall preference averages and subjective mixtures for the class of biseparable preferences,
which was introduced and axiomatized by \cite{GhirardatoMarinacci01}.

For $A \in \Sigma$, define $t\,A\,x :=  t\,I_{A}(s)+x\,I_{A^{c}}(s)$.
Here, $t,x\in\mathbb{R}$, $I_{A}$ is the indicator of event $A$ and $A^{c}$ is the complement of $A$.
Let $U$ be a standard uniform random variable on the probability space associated with the second stage 
and define the set $E_y: =\{y<U<1-y\}$, $y\in[0,1/2]$.
For $t>x$, $v_y:=t\,E_y\,x=t\,I_{E_y}+x\,I_{E_y^c}$ converges to $t$ as $y$ tends to $0$ and to $x$ as $y$ tends to $1/2$.
By continuity as assumed in A3,\footnote{Let $u$ be the midpoint between $t$ and $x$. 
Define $y^*:=\inf \{y|  v_y\succ u \}$. 
The set over which the infimum is taken is not empty because, by A4, $y=0$ is in this set.
As by A3 the set $M=\{v|  u \succ v \}$ is open, $v_{y^*}=\lim_n v_ {y^*+1/n} \in M^c$ cannot be in $M$, 
and thus $v_ {y^*}\succeq u$. 
Suppose 
that $v_{y^*}\succ u$. 
Then, by the definition of $y^*$, we would have that $v_{y^*}\succ u \succeq v_{y^*-1/n}$
and by passing to the limit using A3 arrive at a contradiction. 
Hence, $v_{y^*}\sim u$.} there must then be a $y^*\in (0,1/2)$ such that
$t\succ t\, E_{y^*}\, x \succ x$.
Any event $E$ satisfying this preference ordering for some $t>x$ will in the sequel be called an \textit{essential event}.

We call our binary relation $\succeq$ on $V_{0}$ a \textit{biseparable preference}
if it has a nontrivial representation $\mathcal{U}:V_{0}\rightarrow\mathbb{R}$ for which:
\begin{itemize}
\item[(i)] there exists $\rho:[0,1]\rightarrow[0,1]$ with $\rho$ not identical zero or one
on $(0,1)$ such that, for all $t\succeq x$ and all $A\in\Sigma$,
\begin{equation*}
\mathcal{U}(t\,A\,x)=\phi(t)\rho(P(A))+\phi(x)\left(1-\rho(P(A))\right),
\end{equation*}
where $\phi(t)\equiv \mathcal{U}(t)$ for all $t\in \mathbb{R}$ is increasing;
\item[(ii)] $\mathcal{U}(\mathbb{R})$ is convex.
\end{itemize}



Averages of, say, $t$ and $x$ with respect to a symmetric function $h:\mathbb{R}^2\to \mathbb{R}$ can often be defined as the element $z$ such that $h(t,x)=h(z,z)$, 
with the arithmetic and geometric average corresponding to $h(t,x)=t+x$ and $h(t,x)=t\cdot x$.
Now given $t,x\in \mathbb{R}$, if $t\succeq x$ we say that $y\in \mathbb{R}$ is a \textit{preference average}
of $t$ and $x$ (given an essential event $E$) if $t\succeq y\succeq x$
and
\begin{equation*}
t\,E\,x\sim m_{tEt}\,E\,m_{xEx}\sim m_{tEy}\,E\,m_{yEx}.
\end{equation*}


Next, fix some event $B \in \Sigma $, 
and then construct \textit{state by state} a random variable in $V_{0}$
such that every state $s$ yields the certainty equivalent of the bet $v(s)\,B\,u(s)$ for $v,u \in V_{0}$.
That is, the \textit{state-wise B-mixture} of $v$ and $u$ 
is the random variable $\mathrm{sub}(vBu) \in V_{0}$ formally defined as follows:
For all $s \in S$,
\begin{equation*}
\mathrm{sub}(vBu)(s)\equiv m_{(v(s)Bu(s))}.
\end{equation*}
Provided that $v$ dominates $u$ state-wise, the constructed random variable yields state-wise indifference to bets on the event $B$;
the random variable thus constructed can via CDFs be identified with a subjective mixture of lotteries.

Subsequently, we say that two random variables $\tilde{v},\tilde{u}\in\tilde{V}$
are \emph{comonotonic} if, for every $w\in W$ and every $s,s'\in S$,
\begin{equation*}
(\tilde{v}^w(s')-\tilde{v}^w(s))(\tilde{u}^w(s')-\tilde{u}^w(s))\geq 0.
\end{equation*}
Comonotonic random variables do not provide hedging potential
because their realizations move in tandem without generating offsetting possibilities (\citealp{Schmeidler86,Schmeidler89}, \citealp{Yaari87}).

The next axiom is a weak version of the independence axiom of EU for binary comonotonic random variables, 
discussed in \cite{GhirardatoMarinacci01} (for horse lotteries). 
Extending \cite{ChewKarni94}, it weakens in particular Schmeidler's comonotonic independence axiom, 
and can be traced back to \cite{Gul92} and \cite{Nakamura90}.
\vskip 0.14cm
\textit{A{\sc xiom} A8---Binary Comonotonic Independence:\footnote{In \cite{GhirardatoMarinacci01}, this axiom is referred to as ``Binary Comonotonic Act Independence''. 
As in our setting it is imposed on random variables (rather than acts, i.e.,  horse lotteries), we omit the term ``act''.}  
For every essential $A \in \Sigma$, every $B \in \Sigma$, and for all $v, u, r \in V_{0}$ and $x,x',x'',y,y',y''\in\mathbb{R}$
such that $v=x\,A\,y$, $u= x'\,A\,y'$ and $r=x''\,A\,y''$, if $v, u, r$ are pairwise comonotonic, and
$  x,x'  \geq  x''$
and $y,y' \geq y''$
(or $x'' \geq  x,x' $
and $y'' \geq   y,y' $), then}
\begin{equation*}
v \succeq u \Longrightarrow \mathrm{sub}(vBr)\succeq \mathrm{sub}(uBr).
\end{equation*}

The following proposition is suitably adapted from \cite{GhirardatoMarinacci01}, Theorem~11, and 
\cite{GhirardatoMaccheroniMarinacciSiniscalchi03}, Proposition~1:
\begin{proposition}\label{prop:app}
Suppose that A1-A4 and A8 hold.
Then $\succeq$ is a biseparable preference on $V_{0}$, where $\rho$ is unique and $\phi$ is continuous and unique up to a positive affine transformations.
Furthermore, for $t,x\in \mathbb{R}$ and each essential event $E\in\Sigma$,
$y\in \mathbb{R}$ is a preference average of $t$ and $x$ given $E$
if and only if
\begin{equation*}
\phi(y)=\frac{1}{2}\phi(t)+\frac{1}{2}\phi(x).
\end{equation*}
Hence, preference averages of $t$ and $x$ given $E$ exist for every essential event $E\in\Sigma$,
and they are independent of the choice of $E$ and of the normalization of $\phi(\cdot)$.
\end{proposition}

If $z\in\mathbb{R}$ is a preference average of $t$ and $x$,
then, by the proposition above,
we may define $y$ as a preference average of $t$ and $z$
satisfying $\phi(y)=(3/4)\,\phi(t)+(1/4)\,\phi(x)$.
Using standard preference continuity arguments, we may also define $\alpha:(1-\alpha)$-utility mixtures, $\alpha\in[0,1]$.
\textit{Subjective mixtures of random variables} on the space $\tilde{V}_{0}$ can then be defined point-wise:
Given $\tilde{v},\tilde{u}$ and $\alpha\in[0,1]$, the \textit{subjective mixture}
$\alpha \tilde{v}\oplus\,(1-\alpha)\tilde{u}$ is the random variable $\tilde{r}$
defined by
\begin{equation}\label{subjectivemixture1}
\tilde{r}(w):=\alpha \tilde{v}(w)\oplus\,(1-\alpha)\tilde{u}(w),
\end{equation}
for any $w\in W$, where $\tilde{r}$ satisfies
\begin{equation}
\label{subjectivemixture2}\phi(\tilde{r}(w)) = \alpha \phi (\tilde{v}(w)) + (1-\alpha)\phi(\tilde{u}(w)).
\end{equation}

Next, let us define \textit{subjective additions} of random variables.
We first state the following definition.
\begin{definition}
We call $z$ a preference doubling of $x$ (given an essential event $E$) and write $z=2 \otimes x$
if, in case $x \geq 0$, we have $z \succeq x$ and
$z\,E\,0\sim m_{zEx}\,E\,m_{xE0}$;
or, in case $x < 0$, we have $ x \succeq z$ and
$0\,E\,z\sim m_{0Ex}\,E\,m_{xEz}$.
\end{definition}

\begin{proposition}\label{prop:doubling}
$z=2 \otimes x$ if, and only if,
\begin{equation}
\frac{1}{2}\phi(z)+\frac{1}{2}\phi(0)=\phi(x).
\label{eq:z}
\end{equation}
In particular, $z$ is invariant under positive affine transformations of the function $\phi$,
and for all functions $\phi$ with $\phi(0)=0$, we have $\phi(z) = 2\phi(x)$.
\end{proposition}

Finally, we define point-wise
\begin{equation}\label{subjectivmixture3}
\tilde{v}(w)\oplus \tilde{u}(w) := 2 \otimes (\frac{1}{2} \tilde{v}(w) \oplus \frac{1}{2} \tilde{u}(w)),
\end{equation}
and note that, for any $\phi$ from Proposition~\ref{prop:app} with $\phi(0)=0$,
\begin{equation}
\label{subjectivemixture4}
\phi(\tilde{v}(w)\oplus \tilde{u}(w))=\phi (\tilde{v}(w))+ \phi(\tilde{u}(w)).
\end{equation}

\subsection{New Axioms}
We replace the uncertainty aversion axiom of
\cite{GilboaSchmeidler89} and \cite{MMR06} (Axiom
A6MEU) by the following assumption:
\vskip 0.14cm
\textit{A{\sc xiom} A6---Dual Ambiguity 
Aversion: If $v',u'\in V'_0$ and
$\alpha\in(0,1)$, then
$v'\sim u'$ implies $\a v'\oplus (1-\a) u'\succeq v'$.}
\vskip 0.14cm
Observe that, different from Axiom~A6MEU (as well as Axioms~A7MEU and~A7VP) that applies to mixtures of conditional CDFs,
Axiom~A6 considers subjective mixtures of random variables, introduced in Section~\ref{sec:submixrv}. 
By using ``$\oplus$'' rather than ``$+$'', the utility profiles rather than the random variables themselves are combined;
the DM takes subjective mixtures of random variables that carry no risk (are in $V'_0$) and that she is indifferent to.
A dual ambiguity 
averse DM, then, prefers the ``diversified'', convex combination of the utility profiles of the two random variables ($\a v'\oplus(1-\a) u'$)
to the original non-diversified random variable ($v'$ or $u'$). 
Note that, as a special case when $\phi=\mathrm{id}$, ``$\oplus$'' corresponds to ``+''. 
Hence, one can anticipate that the preference for diversification induced by A6 stems from attitude towards ambiguity rather than from attitude towards wealth. 
In fact, consistent with experiments, our theory does not assume a specific attitude towards wealth such as a globally concave $\phi$. 

Before stating our new Axiom~A7, we briefly discuss ``dual independence'', which inspired it. 
\cite{Yaari87}, in the context of decision under risk, asserts that a preference of $v$ to $u$ induces a preference of $\alpha v + (1-\alpha)r$ to $\alpha u + (1-\alpha)r$, $\alpha \in (0, 1)$, in case $v, r$ and $u, r$ are \textit{pairwise comonotonic} (pc). 
(
This assertion also implies (\textit{ceteris paribus}) a preference of $v+r$ to $u+r$.) 
In particular, \cite{Yaari87} replaces the independence
axiom of EU by the following assumption (involving outcome mixtures of random variables), restricted to decision under risk:
\vskip 0.14cm
\textit{A{\sc xiom} A7D---Dual Independence: 
Let $v, u, r \in V_0$ and assume that $v, r$ and $u, r$ are pc.
Then, for every $\alpha \in (0, 1)$, $v\succeq u\Rightarrow \alpha v + (1-\alpha)r \succeq \alpha u + (1-\alpha)r$.}
\vskip 0.14cm
\noindent To extend this axiom to subjective mixtures of random variables, ``$+$'' in A7D would have to be replaced by ``$\oplus$''.

Against this background, first consider $\tilde{v},\tilde{u},\tilde{r}\in\tilde{V}_0$
and suppose that the DM prefers $\tilde{v}$ to $\tilde{u}$.
Is it natural to require that the DM then also prefers $\tilde{v}\oplus\tilde{r}$ to $\tilde{u}\oplus\tilde{r}$,
or $\a \tilde{v}\oplus(1-\a)\tilde{r}$ to $\a \tilde{u}\oplus(1-\a)\tilde{r}$ with $\a\in(0,1)$,
in general (without comonotonicity imposed)?
If $\tilde{u}$ and $\tilde{r}$ are not comonotonic, then the DM may try to employ $\tilde{r}$
to hedge against adverse realizations of $\tilde{u}$.
As a result, $\tilde{u}\oplus\tilde{r}$ can conceivably be ``better diversified'' than $\tilde{v}\oplus\tilde{r}$
(depending on the joint stochastic nature of $\tilde{u},\tilde{r}$ on the one hand and that of $\tilde{v},\tilde{r}$ on the other),
and the DM may instead prefer $\tilde{u}\oplus\tilde{r}$ to $\tilde{v}\oplus\tilde{r}$.
Hence, consider $\tilde{v},\tilde{u},\tilde{r}\in\tilde{V}_0$, suppose that
$\tilde{v}\succeq \tilde{u}$, as before,
and suppose furthermore that $\tilde{v},\tilde{r}$ and $\tilde{u},\tilde{r}$ are pc.
Should a DM then also prefer $\tilde{v}\oplus\tilde{r}$ to $\tilde{u}\oplus\tilde{r}$?
(Note that this is not implied by Axiom A7D, which requires the random variables to live in the space $V_0$, and besides considers ``$+$'' instead of ``$\oplus$''.)
Even though (subjectively) adding $\tilde{r}$ does, in view of the pc assumption, not induce any discriminatory hedging potential,
it may still impact the ambiguity ``level'', in a discriminatory manner, leading to a preference reversal.

Consider the following example, where for ease of exposition, we assume an affine function $\phi$ so that ``$\oplus$'' and ``$+$'' agree:
\begin{Example}\label{ex:2}
Consider two urns, $A$ and $B$, and
$50$ balls, 25 of which are red and 25 of which are black.
Every urn contains 25 balls.
The exact number of balls per color in each urn is unknown.
Furthermore, consider two urns, $C$ and $D$, and
$50$ balls, 30 of which are red and 20 of which are black.
As for $A$ and $B$, every urn contains 25 balls,
but the exact number of balls per color in each urn is unknown.

Denote by $p_{i}$ the (unknown) probability of drawing a red ball from urn $i, i\in\{A,B,C\}$.
Draw a random number, say $U$, from the set $\{1,\ldots,25\}$ with each number having the same likelihood.
Consider:
\begin{itemize}
\item[(i)] the random variable $\tilde{v}$ that represents a payoff of $\$100$ if $U$ is smaller than or equal to the number of red balls in urn $C$, and 0 otherwise;
\item[(ii)] the random variable $\tilde{u}$ that represents a payoff of $\$100$ if $U$ is smaller than or equal to the number of red balls in urn $A$, and 0 otherwise;
\item[(iii)] the random variable $\tilde{r}$ that represents a payoff of $\$100$ if $U$ is smaller than or equal to the number of red balls in urn $B$, and 0 otherwise.
\end{itemize}
Note that $\tilde{v},\tilde{r}$ and $\tilde{u},\tilde{r}$ are pc.
Typically, $\tilde{v}\succeq\tilde{u}$ because $30>25$.
At the same time, the DM may prefer $\tilde{u}+\tilde{r}$ to $\tilde{v}+\tilde{r}$,
because the former is, loosely speaking, less ambiguous than the latter.
More specifically, the unknown probability of drawing red from $A$ is connected (complementary) to the unknown probability of drawing red from $B$:
with certainty, $p_{A}+p_{B}=1$.
By contrast, the probability of drawing red from $B$ (or $A$) is not connected to the probability of drawing red from $C$.
Mathematically, $\tilde{u}+\tilde{r}$ yields at least $\$100$ with probability $\max\{p_{A},1-p_{A}\}\geq$ 50\%,
and it yields exactly $\$200$ with probability $\min\{p_{A},1-p_{A}\}=1-\max\{p_{A},1-p_{A}\}$.
On the other hand, $\tilde{v}+\tilde{r}$ has potential realizations  $\$0$, $\$100$, and $\$200$
with unknown probabilities, where no non-trivial upper or lower bounds can be given.
$\qquad\triangledown$
\end{Example}

We will assert that, if $\tilde{v},\tilde{u},\tilde{r}\in\tilde{V}_0$,
$\tilde{v}\succeq \tilde{u}$, and $\tilde{v},\tilde{r}$ and $\tilde{u},\tilde{r}$ are pc,
then the implication
$\tilde{v}\oplus\tilde{r}\succeq \tilde{u}\oplus\tilde{r}$
only holds if $\tilde{r}$ carries no ambiguity (i.e., is in $V_{0}$),
hence cannot impact the ambiguity level, in a discriminatory manner.
This motivates to replace the weak certainty independence axiom by the following assumption:
\vskip 0.14cm
{\it A{\sc xiom} A7---Weak Certainty Dual Independence: Let $\tilde{v},\tilde{u}\in
\tilde{V}_0 $ and $r\in V_0$.
Suppose that $\tilde{v},r$ and $\tilde{u},r$ are pc.
Then,}
$\tilde{v}\succeq \tilde{u} \Rightarrow \tilde{v}\oplus r\succeq \tilde{u}\oplus r.
$



\subsection{Main Results}

The following axiomatic characterization generalizes the rank-dependent utility model to the setting of risk \textit{and} ambiguity.
The result may also be viewed as a generalization of the variational preferences model; see Table~\ref{tab:dual}.
As there are no concavity restrictions on the utility function,
and no convexity restrictions on the probability weighting function,
the result encompasses (inverse) S-shaped utility and probability weighting functions used in prospect theory
and supported by empirical evidence.

Let $\psi:[0,1]\rightarrow [0,1]$ be a non-decreasing and continuous
function satisfying $\psi(0)=0$ and $\psi(1)=1$.
We refer to $\psi$ as a probability weighting or distortion function.
For $v\in V$, we define the measure $\nu_{\psi}$ through
\begin{equation*}
\int v \,\mathrm{d}\nu_{\psi}:= \int_{-\infty}^0 (\psi(1-F_v(t))-1)\,\mathrm{d}t +\int_0^\infty \psi(1-F_v(t))\,\mathrm{d}t.
\end{equation*}
One readily verifies that, for $a>0$ and $b\in\mathbb{R}$,
$\int (av+b) \,\mathrm{d}\nu_{\psi}=a\int v \,\mathrm{d}\nu_{\psi}+b$. 
We now state our main result, which provides a representation theorem
characterizing a preference relation satisfying Axioms A1-A8:


\begin{theorem}\label{th:RDU}
\begin{itemize}
\item[$(\alpha)$] A preference relation $\succeq$ satisfies A1-A8 if, and only if,
there exist a non-constant, non-decreasing and continuous function $\psi:[0,1]\rightarrow[0,1]$ with $\psi(0)=0$ and $\psi(1)=1$,
an increasing and continuous function $\phi:\mathbb{R} \rightarrow \mathbb{R}$
and a grounded, convex and lower-semicontinuous function $c:\Delta(W,\Sigma^\prime)\rightarrow[0,\infty]$ such that,
for all $\tilde{v},\tilde{u}\in\tilde{V}_0$,
\begin{align}\begin{split}
\tilde{v}\succeq\tilde{u}\Leftrightarrow&\min_{Q\in\Delta(W,\Sigma^\prime)}\left\{ \mathbb{E}_Q\left[\int \phi(\tilde{v}^\centerdot) \,\mathrm{d}\nu_\psi\right]+c(Q)\right\} 
\\
&\geq \min_{Q\in\Delta(W,\Sigma^\prime)}\left\{ \mathbb{E}_Q\left[\int \phi(\tilde{u}^\centerdot) \,\mathrm{d}\nu_\psi\right]+c(Q)\right\}.
\end{split}\label{eq:rduamb}\end{align}
Furthermore, for each $\phi$ there exists a unique minimal $c_{\min}$ satisfying \eqref{eq:rduamb} given by
\begin{equation*}
c_{\min} (Q) = \sup_{v^\prime\in V^\prime_0}\left\{ m_{\phi(v^\prime)}-\mathbb{E}_Q[\phi(v^\prime)]\right\}.
\end{equation*}
\item[$(\beta)$] There exists a unique extension of $\succeq$ to $\tilde{V}$ satisfying A1-A8 on $\tilde{V}$ and \eqref{eq:rduamb}.
\end{itemize}
\end{theorem}
\noindent(Here, $\tilde{v}^{\centerdot}$ denotes the random variable given by $s\mapsto\tilde{v}^{\centerdot}(s)$.\footnote{In fully explicit form, the numerical representation in \eqref{eq:rduamb} may thus be expressed as 
\begin{align*}
\min_{Q\in\Delta(W,\Sigma^\prime)}\left\{ \int_{W}\int_{S} \phi(\tilde{v}^w(s)) \,\nu_\psi(\mathrm{d}s)\, Q(\mathrm{d}w)+c(Q)\right\}.
\end{align*}
})

A natural question is whether it is possible to restrict the class of \textit{penalty} functions on $\Delta(W,\Sigma^\prime)$ to more specific ones, 
such as a penalty function that only takes the values zero or infinity, and what this would entail behaviorally. 
To this end, we replace (\textit{ceteris paribus}) Axiom A7 by the following stronger (i.e., more restrictive) assumption:

\vskip 0.14cm
{\it A{\sc xiom} A7$^{0}$---Certainty Dual Independence:
Let $\tilde{v},\tilde{u}\in \tilde{V}_0 $ and $r\in V_{0}$.
Suppose that $\tilde{v},r$ and $\tilde{u},r$ are pc.
Then,
$\tilde{v}\succeq \tilde{u}\,\,\, \Leftrightarrow\,\,\,\alpha\tilde{v}\oplus(1-\alpha)r\succeq \alpha\tilde{u}\oplus(1-\alpha)r\ \mbox{ for all}\ \alpha\in(0,1).
$
}
\vskip 0.14cm
If Axioms~A1-A6, A7$^{0}$ and A8 hold, then we
obtain the rank-dependent generalization of the popular \cite{GilboaSchmeidler89} maxmin expected utility representation:

\begin{theorem}\label{th:GS}
\begin{itemize}
\item[(a)] A preference relation $\succeq$ satisfies A1-A6, A7$^{0}$ and A8 if, and only if,
there exist a non-constant, non-decreasing and continuous function $\psi:[0,1]\rightarrow[0,1]$
with $\psi(0)=0$ and $\psi(1)=1$, an increasing and continuous function $\phi:\mathbb{R} \rightarrow \mathbb{R}$
and a non-empty, closed and convex set $\mathcal{Q}\subset \Delta(W,\Sigma')$
such that, for all $\tilde{v},\tilde{u}\in\tilde{V}_{0}$,
\begin{equation}\label{eq:gsamb}
\tilde{v}\succeq \tilde{u}\Leftrightarrow \min_{Q\in\mathcal{Q}}\mathbb{E}_Q\left[ \int \phi(\tilde{v}^{\centerdot}) \,\mathrm{d}\nu_{\psi}\right]
\geq \min_{Q\in\mathcal{Q}}\mathbb{E}_Q\left[ \int \phi(\tilde{u}^{\centerdot})\,\mathrm{d}\nu_{\psi}\right].
\end{equation}
Furthermore, there exists a unique extension of $\succeq$ to $\tilde{V}$
satisfying A1-A6, A7$^{0}$ and A8 on $\tilde{V}$ and \eqref{eq:gsamb}.
\item[(b)] If moreover the numerical representation in \eqref{eq:gsamb} is continuous from below, then
$\mathcal{Q}\subset \Delta_{\sigma}(W,\Sigma')$,
i.e., the minimum may be taken over a convex set of probability measures.
\end{itemize}
\end{theorem}

\subsection{Interpretation}

Define $U$ as the numerical representation in \eqref{eq:rduamb}, i.e.,
\begin{equation*}
U(\tilde{v})=\min_{Q \in \Delta(W,\Sigma')}\left\{\E_Q\left[\int \phi(\tilde{v}^{\centerdot})
\,\mathrm{d}\nu_{\psi}\right]+c(Q)\right\}.
\end{equation*}
The numerical representation may be given the following interpretation.
The function $c$ is non-negative and grounded, i.e.,
for every $Q\in\Delta(W,\Sigma')$, $c(Q)\geq 0$, and there exists at least one
measure $P'\in\Delta(W,\Sigma')$ such that $c(P')=0$.
This measure $P'$ may be considered to be the DM's candidate model (i.e., ``approximation''), selected from the set of all measures
on $(W,\Sigma')$.
If the DM believes that $P'$ is a good (reliable) candidate model,
then she can simply take the ($P'$-)expectation over all evaluations of an objective lottery,
i.e., calculate the rank-dependent $\E'[\int \phi(\tilde{v}^{\centerdot}) \,\mathrm{d}\nu_{\psi}]$,
which would correspond to $c(Q)=\infty$ if $Q\ne P'$, $Q\in\Delta(W,\Sigma')$, in \eqref{eq:rduamb}.
In many situations, however, the DM may not fully trust his candidate model $P'$ and
takes other measures on $(W,\Sigma')$ into account.

One way to proceed would be to assume a worst case approach and consider
the representation
$\min_{Q \in \Delta(W,\Sigma')}\E_Q[\int \phi(\tilde{v}^{\centerdot}) \,\mathrm{d}\nu_{\psi}]=\min_w \int  \phi(\tilde{v}^w) \,\mathrm{d}\nu_{\psi}$,
which corresponds to $c(Q)=0$ for all $Q\in \Delta(W,\Sigma')$, as in Theorem~\ref{th:GS}.
In this case, the DM would consider all measures on $(W,\Sigma')$ equally plausible. 
Alternatively, the DM may consider his candidate model $P'$ to be more plausible than other measures,
but still wants to take other measures into account (non-trivially).
In this case, she would take the minimum over all measures in $\Delta(W,\Sigma')$,
and ``penalize'' every measure $Q$ not equal to $P'$ by a penalty $c(Q)$.
This penalty depends on the degree of plausibility that the DM associates to the measure $Q$.
The function $c$ is therefore often referred to as an ambiguity index.
Such procedures that explicitly account for the fact that the measure $P'$ is only an approximation and may deviate from the true measure
are often referred to as robust approaches.
They are robust against ``malevolent nature''.

In statistics,
engineering and optimal control,
risk management, and robust macroeconomics,
the plausibility of the measure $Q$ is often expressed by the relative
entropy of $Q$ with respect to the approximation $P'$;
see \cite{Csiszar75}, \cite{Ben-Tal85}, \cite{HansenSargent00,
HS01}, \cite{MMR06}, \cite{Strzalecki11a} and \cite{LaevenStadje13, LaevenStadje14}.
In our setting, this would lead to $c(Q)=\theta R(Q|P')$ with
$R(Q|P')=\E_Q\left[\log\left(\frac{\mathrm{d}Q}{\mathrm{d}P'}\right)\right]$ and $\theta$
a non-negative constant.
The relative entropy, also referred to as Kullback-Leibler divergence,
measures the deviation of $Q$ from $P'$ and is zero if and only if $Q\equiv P'$.
Thus, measures that are close to $P'$ are penalized weakly,
while measures that deviate strongly from $P'$
are penalized strongly.
Specifically, for $v'\in V'_0$,\begin{equation*}\min_{Q \in
\Delta(W,\Sigma')}\left\{\E_Q[v']+\theta
R(Q|P')\right\}=-\theta\log\left\{\E'\left[\exp(-v'/\theta)\right]\right\}.\end{equation*}
In general, for $\tilde{v}\in \tilde{V}_{0}$,
\begin{equation*}\min_{Q \in \Delta(W,\Sigma')}\left\{\E_Q\left[\int \phi(\tilde{v}^{\centerdot}) \,\mathrm{d}\nu_{\psi}\right]+c(Q)\right\}=
-\theta\log\left(\E'\left[\exp\left(\int -\phi(\tilde{v}^{\centerdot}) \,\mathrm{d}\nu_{\psi}/\theta\right)\right]\right).
\end{equation*}
Other ways of penalizing ``deviating beliefs'' include
$c(Q)=\theta G(Q|P')$ with $G(Q|P')$
the relative Gini index given by
$\E\left[\left(\frac{\mathrm{d}Q}{\mathrm{d}P'}-1\right)^2\right]$.
Again, $G(Q|P')=0$ if and only if $Q=P'$; $G(Q|P')$ measures how much the ratio of $Q$
and $P'$ deviates from one; see e.g., \cite{MMRT04}.
\cite{MMR06} also propose to weight every state of the world $w$
by a weighting function $h:W\rightarrow \mathbb{R}_+$ satisfying $\int_W h(w)P'(\mathrm{d}w)=1$.
The penalty functions are then given by
$c(Q)=\int_W h(w)\log\left(\frac{\mathrm{d}Q}{\mathrm{d}P'}(w)\right) Q(\mathrm{d}w)$ and
$c(Q)=\int_W h(w)\left(\frac{\mathrm{d}Q}{\mathrm{d}P'}(w)-1\right)^2Q(\mathrm{d}w)$.\footnote{Note that these two penalty functions are not
\textit{probabilistically sophisticated} on $(W,\Sigma')$ unless $h\equiv 1$; see also Section~\ref{sec:ambrisk}.}

\setcounter{equation}{0}

\section{Probabilistic Sophistication and Ambiguity Aversion}\label{sec:ambrisk}

\subsection{Probabilistic Sophistication}

We refer to $P'$ as a \textit{reference measure} on $(W,\Sigma')$ if the DM is indifferent
between random variables that have the same probability distribution under $P'$.
We say that a DM who adopts a reference measure on $(W,\Sigma')$ 
is \textit{probabilistically sophisticated};
see \cite{MachinaSchmeidler92} and \cite{Epstein99}.
Our axioms do not (necessarily) imply the existence of a reference measure on $(W,\Sigma')$.
But
%
in case there is a reference measure $P'$ on $(W,\Sigma')$, we define, for a given $v'\in V_{0}'$,
$F'_{v'}$ by $F'_{ v'}(t)=P'[v'\leq t]$.
With slight abuse of notation we say that $v'\succeq_2 u'$ if, for every $t\in\mathbb{R}$,
$\int _{-\infty}^t F'_{v'}(\tau)\,\mathrm{d}\tau\leq \int_{-\infty}^t F'_{u'}(\tau)\,\mathrm{d}\tau$.
We call $\succeq_2$ \textit{second order stochastic dominance} (SSD) on $V'$ with respect to $P'$;
see \cite{RothschildStiglitz70}. 
More generally, we say that $v'\succeq_{\phi,2} u'$ if, for every $t\in\mathbb{R}$,
$\int _{-\infty}^t F'_{\phi(v')}(\tau)\,\mathrm{d}\tau\leq \int_{-\infty}^t F'_{\phi(u')}(\tau)\,\mathrm{d}\tau$ and call $\succeq_{\phi,2}$ \textit{$\phi$-second order stochastic dominance} ($\phi$-SSD) on $V'$ with respect to $P'$. 
The following proposition shows that the availability of a non-atomic reference measure $P'$ 
is equivalent to requiring
that the DM respects $\phi$-SSD on $V'_0$:\footnote{Recall that $P'$ is non-atomic if the probability space $(W,\Sigma')$ is rich enough to support a uniformly distributed random variable.}$^{,}$\footnote{The preference relation $\succeq$ respects $\phi$-SSD on $V'_0$ if,
for all $v',u'\in V'_0$ with $v'\succeq_{\phi,2} u'$, $v'\succeq u'$.}
\begin{proposition}\label{prop:SSD}
Suppose that a preference relation $\succeq$ satisfies A1-A8.
Then the following statements are equivalent:
\begin{itemize}
\item[(a)] there exists a non-atomic reference measure, say $P'$, on $(W,\Sigma')$.
\item[(b)] $\succeq$ respects $\phi$-SSD on $V'_0$ with respect to a non-atomic $P'$.
\end{itemize}
Furthermore, if $\phi$ is concave,
$\succeq$ admitting a (non-atomic) reference measure $P'$ is equivalent to 
$\succeq$ respecting SSD on $V'_0$.
\end{proposition}


\subsection{Ambiguity Aversion}\label{sec:aa}

Subsequently, we say that $\succeq$ is \textit{more ambiguity averse} than $\succeq^*$ if, 
for all $\tilde{v}\in \tilde{V}_0$ and $v\in V_0$,
\begin{equation*}
\tilde{v}\succeq v\Rightarrow \tilde{v}\succeq^* v.
\end{equation*}
Similar definitions of \textit{comparative ambiguity aversion} can be found,
for instance, in \cite{Epstein99}, \cite{GhirardatoMarinacci02} and
\cite{MMR06}; see also the early \cite{Yaari69}, \cite{Schmeidler89} and \cite{GilboaSchmeidler89}.\footnote{The difference between definitions of uncertainty aversion consists primarily
in the ``factorization'' of ambiguity attitude and risk attitude.
\cite{Schmeidler89} and \cite{GilboaSchmeidler89} adopt the Anscombe-Aumann framework
with objective unambiguous lotteries.
\cite{Epstein99}, by contrast, instead of adopting a two-stage setup and assuming that there exists a
space of objective lotteries, models ambiguity by assuming that
there exists a set of events $\mathcal{A}$ that \emph{every DM}
considers to be unambiguous.
Then he defines comparative ambiguity aversion
through the random variables that are measurable with respect to $\mathcal{A}$.
The model-free factorization approach of \cite{GhirardatoMarinacci02}
in principle encompasses both approaches to modeling ambiguity.}
Our notion of more ambiguity averse agrees with the comparative ambiguity aversion concept
in \cite{GhirardatoMarinacci02} and \cite{MMR06}.

The following result characterizes comparative ambiguity aversion in the setting of our main representation result:
\begin{proposition}\label{prop:caa} 
Consider two preference relations, $\succeq$ and $\succeq^*$, 
induced by assuming Axioms~A1-A8. 
Then,
$\succeq$ is more ambiguity averse than $\succeq^*$ if, and only if, $\succeq^*$ and $\succeq$ may be identified with $(\phi^*,\psi^*,c^*)$ and $(\phi,\psi,c)$ such that $\phi^* = \phi$,
$\psi^* = \psi$ and $c^*\geq c$.
\end{proposition}

We note that in the absence of a probability weighting function,  Proposition~\ref{prop:caa} holds similarly in the primal framework of \cite{MMR06}.

In \cite{Epstein99}, \cite{GhirardatoMarinacci02} and \cite{MMR06} a DM is considered to
be \textit{ambiguity averse} if and only if she is more ambiguity averse than
an \textit{ambiguity neutral} DM.
While in \cite{GhirardatoMarinacci02} and
\cite{MMR06} ambiguity neutrality is equivalent to having SEU
preferences, \cite{Epstein99} identifies ambiguity neutrality
with probabilistic sophistication.
\cite{GhirardatoMarinacci02}, however, argue that in full generality (unless the probability space
is rich enough) probabilistically sophisticated behavior may still
include behavior that can be considered to be ambiguity averse.\footnote{
For instance, if $W$ has only finitely many elements, identifying
ambiguity neutrality with probabilistic sophistication would
imply that a DM with a numerical representation of the form 
$U(\tilde{v})=\min_{Q\in \Delta(W,\Sigma')}\{\E_Q[\int \phi(\tilde{v})\,\mathrm{d}\nu_{\psi}]\}=\inf_w \int
\phi(\tilde{v}^w) \,\mathrm{d}\nu_{\psi}$ would be ambiguity neutral,
at least, if $P'$ does not exclude any $w\in W$.
(That is, $P'[w]>0$ for all $w\in W$.)
This seems counterintuitive in our setting, since the ``worst
ambiguity case'' possible is assumed.
A worst case DM is also
probabilistically sophisticated if $W$ is a subset of $\mathbb{R}^d$
and $P'\sim Leb.$ $W\subset \mathbb{R}^d$ is typically satisfied in
a Bayesian framework.
\cite{Strzalecki11b}, however, proves that in the specific framework of \cite{MMR06},
ambiguity neutrality in the sense of \cite{Epstein99},
with non-trivial no-ambiguity sets, implies that the DM has
preferences given by SEU.
\cite{Marinacci02} had proven the same
result under MEU. 
}
Consequently, in our setting, instead of identifying ambiguity
neutrality $(\succeq^{\mathrm{AN}})$ with probabilistic sophistication,
it seems more natural to define $\succeq^{\mathrm{AN}}$ via a numerical representation that induces
computing a plain expectation on the space $W$ with respect
to some measure $P'$.
In other words, we consider a DM to be
\emph{ambiguity neutral} if there exist a measure $P'$, a utility function $\phi$ and a
probability weighting function $\psi$ such that, 
for all $\tilde{v},\tilde{u}\in\tilde{V}_{0}$,
\begin{equation*}
\tilde{v}\succeq \tilde{u} \Leftrightarrow
\mathbb{E}_{P'}\left[\int \phi(\tilde{v}^{\centerdot}) \,\mathrm{d}\nu_{\psi}\right]\geq \mathbb{E}_{P'}\left[\int \phi(\tilde{u}^{\centerdot}) \,\mathrm{d}\nu_{\psi}\right].
\end{equation*}
Next, we say that a DM with a preference relation $\succeq$ is 
\emph{ambiguity averse} if there exists an ambiguity
neutral DM with a preference relation $\succeq^{\mathrm{AN}}$ such that
$\succeq$ is more ambiguity averse than $\succeq^{\mathrm{AN}}$.\\
\begin{proposition}\label{prop:aa}
If $\succeq$ satisfies A1-A8, then $\succeq$ is ambiguity averse.
\end{proposition}

\setcounter{equation}{0}

\section{Motivation and Applications}\label{sec:motivation}

Besides the Allais and Ellsberg paradoxes, which can jointly be rationalized by our theory,
and the disentanglement of attitudes towards ambiguity, risk, and wealth that our theory permits,
we offer three additional motivations for our rank-dependent theory for decision under risk and ambiguity.

\subsection{Robust Risk Management}\label{sec:robustrm}

Our theory provides a decision-theoretic foundation for robust tail risk measures,
unifying tail risk measures as described e.g., in \cite{FollmerSchied16}, Section~4.6, 
and models with robustness of e.g., \cite{HS95, HS01, HS07} type  
or general $f$-divergences.
Indeed, our theory encompasses tail risk measures when probabilities are given, but also accounts for ambiguity.

Tail risk measures are extensively used by financial institutions and regulators
to control and manage risks and to determine adequate capital reserves.
The industry-standard tail risk measures 
are Value-at-Risk ($\text{VaR}_{\lambda}$) and Expected Shortfall ($\text{ES}_{\lambda}$).\
Value-at-Risk is the amount of capital needed to guarantee that
with a certain probability
no shortfall will be suffered over a pre-specified time horizon.\footnote{That is, it is the quantile function of the respective probability distribution; see, e.g., \cite{DP97} and \cite{J97} for a detailed discussion.}
Expected Shortfall measures the expected loss beyond the Value-at-Risk.
More formally, for a random variable $v$ with a given probabilistic model $P$,
$\text{VaR}_{\lambda}(v):=\inf\{t\in\mathbb{R}|P[-v\leq t]\geq 1-\lambda\}$,
$\lambda\in(0,1)$,
and
$\text{ES}_{\lambda}(v):=\frac{1}{\lambda}\int_{0}^{\lambda}\text{VaR}_{\gamma}(v)\,\mathrm{d}\gamma$, $\lambda\in(0,1]$.

As is well-known (see, e.g., \citealp{FollmerSchied16}, Section~4.6),
the numerical representation of Yaari's DT can be expressed as
\begin{align}
	U(v) &= \int_{-\infty}^{0}\left(\psi(P[v>t])-1\right)\,\mathrm{d}t+\int_{0}^{\infty}\psi(P[v>t])\,\mathrm{d}t 
    = \int_{0}^{1} \text{VaR}_{\gamma}(v)\,\mathrm{d}\psi(1-\gamma),
\end{align}
provided either $\mathrm{VaR}_{\lambda}$ or $\psi$ is continuous.
Thus, DT corresponds to ``weighted VaR'' risk measures,
where quantiles are weighted according to the probability weighting function $\psi$.
Expected Shortfall corresponds to (minus) weighted VaR with probability weighting function given by
$\psi(\gamma) = \frac{1}{\lambda}\max( \gamma-(1-\lambda),0)$.\footnote{Other examples of weighted VaR risk measures include the inter-quartile range, the absolute deviation, Gini-related risk measures, Range Value-at-Risk, and dual-power risk measures; see \cite{D94}, \cite{FollmerSchied16}, \cite{EL21} and the references therein. 
}
Clearly, weighted VaR risk measures are measures of ``tail risk''.

Weighted VaR and related risk measures have been widely considered in economics, finance, insurance, operations research and statistics.\footnote{See, e.g., \cite{RV03}, \cite{C06}, 
	\cite{CM09}, \citet{BS01}, 
	\cite{KPH13},
	\cite{ELS20}.}
The primal theories for decision under risk and ambiguity of Table~\ref{tab:dual}, however,
do not provide a foundation for the use of weighted VaR risk measures under a known probabilistic model:
a DM who complies with one of the primal theories adopts the EU model under risk---evaluating the expected utility of the risk or computing the corresponding certainty equivalent or the indifference price.
None of these formally correspond to using a weighted VaR risk measure.

In recent years---in particular after the failure of risk management systems during the global financial crisis---,
it has been argued that by assuming a known probabilistic model,
risk measures such as weighted VaR
fail to take into account that the adopted probabilistic model
may be misspecified: incorrect and only an approximation.
Consequently, there have been calls to seek for risk measures
that account for ambiguity and lead to more robust risk management,
by ensuring that policy performs well under a wide range of potential probabilistic models;
see, e.g., \cite{H14}. 
Such robust risk measures arise in much generality in our representation \eqref{eq:numrep}, taking $\phi$ to be affine. 

For example, robust VaR risk measures with robustness of max-min type, Hansen-Sargent type (i.e., with Kullback-Leibler divergence/relative entropy) or with general $f$-divergences occur as special cases of our numerical representation, as follows:
\begin{align*}
	U(\tilde{v})=\min_{Q \in \Delta}
	\left\{\mathbb{E}_Q\left[\mathrm{VaR}_{\lambda}\left(\tilde{v}^{\centerdot}\right)\right]+c(Q)\right\}.
\end{align*}
Risk measures in this spirit have recently been used e.g., by \cite{GOO03}, \cite{Embrechts13}, \cite{L22}, \cite{CLM22} and \cite{PWW22}.
We provide 
their decision-theoretic foundation.
Indeed, contrary to the primal theories, our new theory for decision under risk and ambiguity
provides a decision-theoretic, axiomatic foundation for such an approach:
it unifies 
weighted VaR risk measures
and the robust approach that explicitly takes ambiguity into account.


\subsection{
	Portfolio Choice with Mean Risk Models}

An important, related motivation for our theory is that it provides a decision-theoretic foundation
for the use of \textit{mean risk} performance criteria in optimal portfolio choice.
In general, mean risk models are not compatible with the primal theories of Table~\ref{tab:dual}, but have seen much interest in theory and applications.

More specifically, letting $\tilde{v}$ denote the uncertain return of a portfolio at time $T$ where $T>0$ is a fixed time horizon,
\citet{RV03}, \cite{W14}, \citet{HX15} and \cite{CLM22}, among many others, 
consider performance criteria of the form
\begin{align}\label{eq:maxmr}
\mathbb{E}_{P}[\tilde{v}]-\rho(\tilde{v})
=\mathbb{E}_{P}\left[\tilde{v}\right]-\min_{Q \in \Delta}
\left\{\mathbb{E}_{Q}\left[\int\tilde{v}^{\centerdot}\,\mathrm{d}\nu_{\psi}\right]+c(Q)\right\},
\end{align}
where $\phi=\text{id}$ and, in the case of mean-(robust) Expected Shortfall, $\psi(\gamma) = \frac{1}{\lambda}\max( \gamma-(1-\lambda),0)$.
In particular, $\rho(\tilde{v})$ is a coherent risk measure in \citet{RV03}
and a robustified convex risk measure in \cite{W14},
both with special emphasis given to (robustified) Expected Shortfall,
while it is a weighted VaR risk measure in \citet{HX15}
and a robustified weighted VaR risk measure in \cite{CLM22}.
These approaches, and, importantly, the diverse optimal investment behavior induced by them, 
are all compatible 
with our theory.
That is, a DM maximizing \eqref{eq:maxmr} can be identified with a DM choosing her optimal portfolio under preferences given by \eqref{eq:numrep}.
Thus, again, our representation provides the decision-theoretic foundation.

Weighted VaR risk measures, with convex $\psi$,
are special cases of coherent (and convex) risk measures,
which have attracted considerable attention in the literature.\footnote{See e.g.,
\citet{ADEH99}, \citet{PCh01}, \cite{FS02}, \cite{FRG02}, 
\citet{Dana05}, \citet{RS06a}, \citet{LNS07}, 
\citet{LaevenStadje13, LaevenStadje14}, \cite{FollmerSchied16} and the references therein.
Specifically, modulo a change of sign, weighted VaR risk measures with convex $\psi$ are coherent risk measures that are law-invariant and additive for comonotonic random variables;
see \cite{K01}.
}
Via their dual representation, 
coherent risk measures can be given the interpretation of accounting for ambiguity.
More specifically, a coherent risk measure can be identified with a DM
who computes the expectation of a random variable with respect to a probabilistic model,
but is uncertain about the probabilistic model, and therefore considers a robust, worst-case expectation
over a collection of probabilistic models.
That is, 
the DM is \textit{risk} neutral but averse to \textit{ambiguity}.
For law-invariant coherent risk measures, such as Expected Shortfall,
this interpretation, while formally correct, may however be unnatural:
when a DM calculates Expected Shortfall,
the probability distribution of the random variable is often assumed to be given,
hence the DM faces decision under risk rather than under ambiguity, and displays risk aversion.
The different possible interpretations of coherent risk measures
arise from the fact that coherent risk measures operate in a one-stage setting, 
where attitudes towards risk are not distinguished from attitudes towards ambiguity.
Our decision theory encompasses coherent and convex risk measures as special cases
while maintaining a clear separation between attitude towards risk and attitude towards ambiguity.

\subsection{Preference for Randomization vs.\ Preference for Diversification}

An important aspect of our theory is that
one of its key axioms entails a preference for \emph{diversification},
while the corresponding axiom in the primal theories is based on 
\emph{randomization}.
In particular, our Axiom A6 stipulates that if a DM is indifferent between two portfolios $A$ and $B$ with unknown probability distributions,
then holding a subjective mixture of $\alpha\%$ of $A$ and $(100-\alpha)\%$ of $B$ is preferred to holding only $A$ or only $B$.
The primal theories, instead, assert that the DM prefers a probabilistic mixture where she receives
portfolio $A$ with $\alpha\%$ probability and portfolio $B$ with $(100-\alpha)\%$ probability
to having for sure only $A$ or only $B$.

A preference for diversification is a trait of behavior that is widespread in financial, insurance, and other real-life situations,
often summarized as ``Don't put all your eggs in one basket''.
It can be traced back to \cite{Bernoulli38}, is a main principle in economics, operations research, and statistics,
and is at the basis of classical optimal investment and asset pricing theories such as Mean-Variance optimization and the Capital Asset Pricing Model; see, e.g., \cite{BE21} for a detailed discussion.
These authors also show that the value of diversification may reduce under ambiguity aversion.
That is, existing definitions of ambiguity aversion, and the primal theories of Table~\ref{tab:dual}, do not necessarily induce a preference for  diversification.
Instead, our Axiom A6 directly implies a diversification preference.

\begin{appendix}

\setcounter{equation}{0}

\renewcommand{\thesection}{\Roman{section}}

\section{Appendix: Proofs} 


\textit{Proof of Proposition~\ref{prop:doubling}.}
\begin{proofwithoutproofatstart}
Assume without loss of generality that $x \geq 0$.
Then,
\begin{align*}
\mathcal{U}\left(m_{zEx}\,E\,m_{xE0}\right)&=\phi(m_{zEx}) \rho(P(E)) + \phi(m_{xE0})\left(1-\rho(P(E))\right) \\
&= \phi(z) (\rho(P(E)))^{2} + 2 \phi(x)\left(1-\rho(P(E))\right)\rho(P(E)) \\
&\ \ \ + \phi(0)\left(1-\rho(P(E))\right)^2 \\
&= \phi(z)\rho(P(E)) + \phi(0)\left(1-\rho(P(E))\right),
\end{align*}
where the last equality holds if, and only if, \eqref{eq:z} holds.
\end{proofwithoutproofatstart}

\textit{Proof of Theorem~\ref{th:RDU}.}
For showing ``$\Rightarrow$'', we let the function $\phi$ be defined through Proposition~\ref{prop:app}, 
whereas for showing ``$\Leftarrow$'', the function $\phi$ is given by the thesis of Theorem~\ref{th:RDU}, in particular, by representation~\eqref{eq:rduamb}.  
By adding or subtracting a constant, we may in both cases assume that $\phi(0)=0$. 
By re-scaling $\phi$ and $c$ if necessary, we may additionally assume that $\phi(1)=1$.
Denote by $\phi^{-1}$ the (increasing  and continuous) inverse of the (increasing and continuous) function $\phi$. 
Since $\phi$ is increasing and continuous, the image of $\phi$, 
henceforth denoted by $\Imp\subset\mathbb{R}$, is an interval. 
Let 
$$\tilde{V}_0^\phi:=\{\tilde{v}\in \tilde{V}_0|\tilde{v}\mbox{ only takes values in } \Imp\},$$
and define $V_0^\phi$ and $V_0^{\prime,\phi}$ similarly. 
Clearly, all these spaces are convex
as $\Imp$ is a convex set. 

Next, for $\tilde{v},\tilde{u}\in \tilde{V}_0^\phi$, we define:\footnote{The auxiliary notation $\succeq^*$ is not to be confused with the preference relation of the less ambiguity averse DM in Section~\ref{sec:aa}.} 
\begin{equation}
\tilde{v} \succeq^* \tilde{u} \text{ if, and only if, } \phi^{-1}(\tilde{v}) \succeq \phi^{-1}(\tilde{u}).
\label{eq:auxiliary}
\end{equation}
We state the following lemma:
\begin{lemma}\label{starpreference}
$\succeq^*$ satisfying A1-A7 with ``$\oplus$'' (in A6-A7) replaced by ``$+$'' and all axioms restricted to random variables taking values in $\Imp$ (i.e., with
$\tilde{V}_0,V_0,V_0^{\prime},m\in\mathbb{R}$ replaced by $\tilde{V}^\phi_0,V^\phi_0,V_0^{\prime,\phi},m\in\Imp$) 
is equivalent to $\succeq$ satisfying A1-A7.
\end{lemma}
\textit{Proof.}
\begin{proofwithoutproofatstart}
For A1-A5, the equivalence is straightforward to see, noting for A3 that by continuity of $\phi$ and $\phi^{-1}$, $\phi(A)$ and $\phi^{-1}(B)$ are open sets if $A\subset\mathbb{R}$ and $B\subset\Imp$ are open. 
Let us show the equivalence for A6 and A7.\\
A6: $\succeq$ ``$\Rightarrow$'' $\succeq^{*}$. 
Assume that A6 holds for $\succeq$ and let us show that A6 also holds for $\succeq^*$. 
By definition of $\succeq^{*}$,
\begin{align*}
&v^\prime \sim^{*} u^\prime \text{ if, and only if, } \phi^{-1}(v^\prime) \sim \phi^{-1}(u^\prime).
\end{align*}
From A6 for $\succeq$, we then have that
$\tilde{r} := \alpha \phi^{-1}(v^\prime) \oplus (1-\alpha)\phi^{-1}(u^\prime) \succeq \phi^{-1}(v^\prime)$. 
Hence, by definition \eqref{subjectivemixture1}--\eqref{subjectivemixture2},
$$ \alpha v^\prime + (1-\alpha)u^\prime=\phi(\tilde{r})\succeq^{*}  \phi(\phi^{-1}(v^\prime))=v^\prime.$$
A7: $\succeq$ ``$\Rightarrow$'' $\succeq^{*}$. 
Recall \eqref{eq:auxiliary}.
Now $r \in V_0^{\phi}$ being comonotonic to $\tilde{v}$ and $\tilde{u}$ is equivalent to $\phi^{-1}(r)$ being comonotonic to $\phi^{-1}(\tilde{v})$ and $\phi^{-1}(\tilde{u})$. 
Hence, \eqref{eq:auxiliary} and A7 entail that
$$\tilde{r}:= \phi^{-1}(\tilde{v}) \oplus \phi^{-1} (r) \succeq \phi^{-1}(\tilde{u}) \oplus \phi^{-1}(r) =: \tilde{\tilde{r}}.$$
Thus, by definition \eqref{subjectivmixture3}--\eqref{subjectivemixture4},
$$\tilde{v}+r =  \phi(\tilde{r})\succeq^* \phi(\tilde{\tilde{r}})=\tilde{u}+r.$$
A6: $\succeq$ ``$\Leftarrow$'' $\succeq^{*}$.  
From $v^\prime \sim u^\prime$ it follows that $\phi (v^\prime) \sim^* \phi(u^\prime)$ and therefore $\alpha \phi(v^\prime) + (1-\alpha) \phi(u^\prime) \succeq^* \phi(v^\prime)$. 
Thus, 
$$\alpha v^\prime \oplus (1-\alpha) u^\prime=\phi^{-1}(\alpha \phi(v^\prime) + (1-\alpha) \phi(u^\prime) ) \succeq  \phi^{-1}(\phi(v^\prime))=v^\prime.$$
A7: $\succeq$ ``$\Leftarrow$'' $\succeq^{*}$.
Clearly, $\phi(r)$ is comonotonic to $\phi(\tilde{v})$ and $\phi(\tilde{u})$. 
It follows by A7 for $\succeq^*$ that $\phi(\tilde{v}) + \phi(r) \succeq^* \phi(\tilde{u}) + \phi(r)$. 
Then we can conclude that
$$\tilde{v} \oplus r =\phi^{-1}(\phi(\tilde{v})+\phi(r))\succeq \phi^{-1}(\phi(\tilde{u})+\phi(r))=\tilde{u} \oplus r.$$ 
\end{proofwithoutproofatstart}

Upon exploiting Lemma~\ref{starpreference} jointly with Proposition~\ref{prop:app}, we can prove Theorem~\ref{th:RDU} by establishing the following result:
\begin{theorem}
\label{1}
\begin{itemize}
\item[($\alpha)$]
A preference relation $\succeq^*$ satisfies A1-A7 with ``$\oplus$'' (in A6-A7) replaced by ``$+$'' and for random variables taking values in $\Imp$ (i.e., with $\tilde{V}_0,V_0,V_0^{\prime},m\in\mathbb{R}$ replaced by $\tilde{V}^\phi_0,V^\phi_0,V_0^{\prime,\phi},m\in\Imp$) if, and only if,
there exist a non-constant, non-decreasing and continuous function $\psi:[0,1]\rightarrow [0,1]$ with
$\psi(0)=0$ and $\psi(1)=1$ and a grounded, convex and 
lower-semicontinuous function $c:\Delta(W,\Sigma')\rightarrow [0,\infty]$	
such that, for all $\tilde{v},\tilde{u}\in\tilde{V}^\phi_0$,
\begin{align}\label{dual2}\begin{split} 
\tilde{v}\succeq^*
\tilde{u}\Leftrightarrow \min_{Q \in \Delta(W,\Sigma')}\left\{\E_Q\left[\int
\tilde{v}^{\centerdot}\,\mathrm{d}\nu_{\psi}\right]+c(Q)\right\}\qquad
\\\geq \min_{Q \in
	\Delta(W,\Sigma')}\left\{\E_Q\left[\int \tilde{u}^{\centerdot}\,\mathrm{d}\nu_{\psi}\right]+c(Q)\right\}.
	\end{split}\end{align}
	Furthermore, there exists a unique 
	minimal $c_{\mathrm{min}}$ satisfying \eqref{dual2} given by
	\begin{equation*}
c_{\mathrm{min}}(Q)=\sup_{v'\in V^{\prime,\phi}_0}\left\{m_{v'}-\E_Q[v']\right\},
\end{equation*}
where $m_{v'}$ is the certainty equivalent of $v'$ under $\succeq^*$.
\item[($\beta$)] There exists a unique extension of $\succeq^*$ to $\tilde{V}^{\phi}$
satisfying A1-A7 on $\tilde{V}^{\phi}$ and \eqref{dual2}, with ``$\oplus$'' in A6-A7 replaced by ``$+$''. 
\end{itemize}
\end{theorem}

\textit{Proof of Theorem~\ref{1}.} 
\begin{proofwithoutproofatstart}
Using for A7 that the integral w.r.t.~$\mathrm{d}\nu_{\psi}$ in \eqref{dual2} is additive for comonotonic random variables, the only property that is
not straightforward to verify in the ``if'' part of Theorem~\ref{1}$(\alpha)$
is the continuity property (Axiom A3).
Let $U$ be the numerical representation in (\ref{dual2}).
This implies that, for all $v\in V_0^{\phi}$,
\begin{equation*}U(v)=\int_{-\infty}^0(\psi(1-F_{v}(t))-1)\,\mathrm{d}t+\int_0^\infty\psi(1-F_v(t))\,\mathrm{d}t.\end{equation*}
The first part of A3 would follow if we could show that $U$ is continuous with respect to weak convergence of uniformly bounded sequences.
Thus, suppose that $v_n$ is a uniformly bounded sequence in $V^\phi_0$, and $v_n\rightarrow v$,
in distribution.
Then, by definition, $F_{v_n}$ converges to $F_v$ at all continuity points of $F_v$.
Because $F_v$ and the $F_{v_n}$'s are non-decreasing functions,
they are continuous, Lebesgue almost everywhere.
But this implies that $F_{v_n}$ converges to $F_v$,
Lebesgue almost everywhere.
Furthermore, because $v_n$ is uniformly bounded by a constant, say $M$,
$F_{v_n}(t)\in \{0,1\}$ for $t\notin [-M,M]$.
In view of the point-wise convergence of $F_{v_n}$ to $F_v$, Lebesgue almost everywhere,
this implies that $F_{v}(t)\in \{0,1\}$ for $t\notin [-M,M]$, as well.
Finally, because $\psi$ is a continuous function, it is bounded on $[0,1]$.
Hence,
\begin{align*}
\lim_n U(v_n)=&\lim_n \int_{-\infty}^0(\psi(1-F_{v_n}(t))-1)\,\mathrm{d}t+\int_0^\infty\psi(1-F_{v_n}(t))\,\mathrm{d}t\nonumber\\
=&\lim_n \int_{-M}^0(\psi(1-F_{v_n}(t))-1)\,\mathrm{d}t+\int_0^M\psi(1-F_{v_n}(t))\,\mathrm{d}t\nonumber\\
=&\int_{-M}^0(\psi(1-F_{v}(t))-1)\,\mathrm{d}t+\int_0^M\psi(1-F_{v}(t))\,\mathrm{d}t
=
U(v),
\end{align*}
as desired.
Proving the second part of Axiom A3 is straightforward and is therefore omitted,
as is the verification of Axioms A1-A2 and A4-A7.

The proof of the ``only if'' part of Theorem \ref{1}$(\alpha)$ consists of the following four steps:
\begin{itemize}
\item[1.] We first show that $\succeq^*$ has a numerical representation $U$ on $\tilde{V}^\phi_0$ 
satisfying certain properties.
\item[2.] Next, we prove that, for all $v\in V^\phi_{0}$, 
$U(v)=\int v \,\mathrm{d}\nu_{\psi}.
$
\item[3.]
Then, we show that, for all $v'\in V^{\prime,\phi}_0$, 
\begin{equation}U(v')=\min_{Q\in \Delta(W,\Sigma')}  \left\{\E_Q[v']+c(Q)\right\}.\label{Urep}\end{equation}
\item[4.] Finally, we derive from Steps 2 and 3 that (\ref{dual2}) holds on $\tilde{V}^\phi_0$.
\end{itemize}

Before proceeding to Step 1, we state the following preliminary lemmata, assuming Axioms A1-A7 hold:

\begin{lemma}
\label{como2}
Let $\tilde{v},\tilde{u},\tilde{v}+m,\tilde{u}+m\in \tilde{V}^\phi_0$ and $m\in \mathbb{R}$.
If $\tilde{v}\succ^*\tilde{u}$ and $\tilde{v},\tilde{u}$ are pc,
then $\tilde{v}+m\succ^*\tilde{u}+m$.
\end{lemma}
\textit{Proof.}
\begin{proofwithoutproofatstart}
By A7, $\tilde{v}+m\succeq^* \tilde{u}+m$.
Suppose that $\tilde{v}+m\sim^* \tilde{u}+m$ would hold.
Then, again by A7, $\tilde{v}=\tilde{v}+m-m \sim^* \tilde{u}+m-m=\tilde{u}$, which is a contradiction.
\end{proofwithoutproofatstart}
\begin{lemma}
\label{certequ}
For every $\tilde{v}\in\tilde{V}_0^\phi$ there exists a certainty equivalent $m^*_{\tilde{v}}\in\Imp$
such that $\tilde{v}\sim^* m^*_{\tilde{v}}$.
\end{lemma}
\textit{Proof.}
\begin{proofwithoutproofatstart}
Suppose, by contradiction, that the lemma does not hold.
Then, the two sets $\{m\in\Imp|m\succ^* \tilde{v}\}$
and $\{m\in\Imp| \tilde{v}\succ^* m\}$
are disjoint, open (by A3) and
their union is equal to $\Imp$.
Furthermore, the two sets are non-empty since, for example, $\max\tilde{v}\in  \Imp$ and $\Imp$ is an open set (as the image of an open set of a strictly increasing function), so that for an $\epsilon>0$ small enough $\max\tilde{v}+\epsilon\in \Imp$. 
Hence, by A4,
$$m:=\max\tilde{v}+\epsilon\succ^* \max\tilde{v}\succeq^* \tilde{v},$$
so that $\{ m\in\Imp| m\succ^* \tilde{v} \}$ is indeed non-empty. 
Similarly, the set $\{m\in\Imp|\tilde{v}\succ^* m\}$ can be seen to be non-empty. 
Because the union of two non-empty disjoint open sets cannot be equal to an open interval (namely $\Imp$), this leads to a contradiction.
Thus, there exists $m^*_{\tilde{v}}\in \mathbb{R}$ such that $\tilde{v}\sim^* m^*_{\tilde{v}}$. 
Since $\min\tilde{v}\leq \tilde{v}\leq \max\tilde{v}$ and $\min\tilde{v},\max\tilde{v}\in \Imp$ (which is an interval), 
it follows that $m^*_{\tilde{v}}\in\Imp$.
\end{proofwithoutproofatstart}

\smallskip
\noindent\textbf{Step 1:}\\
We prove first that $\succeq^*$ has a numerical representation $U:\tilde{V}^\phi_0\rightarrow \mathbb{R}$,
i.e., for all $\tilde{v},\tilde{u}\in\tilde{V}^\phi_{0}$,
\begin{equation*}
\tilde{v}\succeq^* \tilde{u}\Leftrightarrow U(\tilde{v})\geq U(\tilde{u}).
\end{equation*}
We will further show that $U$ 
satisfies the following properties:
\begin{itemize}
\item[(i)] {\it Conditional Law Invariance:} $U(\tilde{v})$ depends only on $F_{\tilde{v}}$.
\item[(ii)] {\it Continuity:} Suppose that $v_n$ is a uniformly bounded sequence in $V^\phi_0$
converging in distribution to $v\in V^\phi_0$,
then $\lim_n U(v_n)= U(v)$.
\item[(iii)] {\it Certainty First-Order Stochastic Dominance:}
For all $v,u\in V^\phi_0$: If $F_{v}(t)\leq F_{u}(t)$ for every $t\in\mathbb{R}$, then $U(v)\geq U(u)$.
\item[(iv)] {\it Monotonicity:} For all $\tilde{v},\tilde{u}\in \tilde{V}^\phi_0$: If $U(\tilde{u}^{w})\leq U(\tilde{v}^{w})$
for every $w\in W$, then $U(\tilde{u})\leq U(\tilde{v})$.
\item[(v)] {\it Certainty Comonotonic Additivity:} Let $\tilde{v}\in\tilde{V}^\phi_{0}$ and $r\in V^\phi_0$.
Suppose that $\tilde{v},r$ are pc and $\tilde{v}+r\in\tilde{V}^\phi_{0}$.
Then $U(\tilde{v}+r)=U(\tilde{v})+U(r)$.
\item[(vi)] {\it Certainty Positive Homogeneity:} For all $v\in V^\phi_0$ and $\lambda\geq 0$,
$U(\lambda v)=\lambda U(v)$ provided that $\lambda v\in V^\phi_0$.
\item[(vii)] {\it Translation Invariance:} For all $\tilde{v}\in\tilde{V}^\phi_0$ and $m\in\mathbb{R}$ such that $\tilde{v}+m\in\tilde{V}^\phi_0$,
$U(\tilde{v}+m)=U(\tilde{v})+m$.
\item[(viii)] {\it Ambiguity 
Concavity}:
If $v',u'\in V^{\prime,\phi}_0$ and $\a\in(0,1)$, then
$U(\a v'+(1-\a)u')\geq \a U(v')+(1-\a)U(u')$.
\end{itemize}

Assume Axioms A1-A7 hold.
For $\tilde{v}\in \tilde{V}^\phi_0$, set $U(\tilde{v})=m^*_{\tilde{v}}$. 
By Lemma~\ref{certequ}, $U$ is well-defined and, by A4,
$m^*_{\tilde{v}}$ is unique.
Note that with this definition,
for all $m\in \Imp$, $U(m)=m$.
Furthermore, it follows from the strict
monotonicity (A4) that
$U(\tilde{v})>U(\tilde{u})$ 
if and only if $\tilde{v}\succ^{*} \tilde{u}$.
Thus, $U$ is a numerical representation of $\succeq^*$.

Next, let us show that $U$ satisfies properties (i)-(viii).
Properties (i)-(iv) with all random variables involved only taking values in $ \Imp$ follow directly from Axioms A1-A5 and the fact that $U$ is a
numerical representation of $\succeq^*$.

To see property (vii), note that Axiom A7 implies that if
$\tilde{v}, \tilde{u} ,
\tilde{v}+m, \tilde{u}+m\in \tilde{V}_0$, we have
\begin{equation}\label{cashinvariance} \tilde{v}\sim^* \tilde{u} \Leftrightarrow
\tilde{v}+m\sim^* \tilde{u}+m.\end{equation}
We claim that this implies that $U$ is translation invariant.
This can be seen as follows. 
If $\tilde{v}\sim^* m^*_{\tilde{v}}$, then, for all $m\in\mathbb{R}$ such that $\tilde{v}+m\in \tilde{V}_0^\phi$, 
$$\min \tilde{v}+m\leq m^*_{\tilde{v}} +m\leq \max \tilde{v}+m.$$
Since $ \tilde{v}+m\in \tilde{V}_0^\phi$, the left-hand side and the right-hand side are both in $\Imp$. 
As $\Imp$ is a (convex) interval, it follows that $m^*_{\tilde{v}} +m\in \Imp$.
Hence, by \eqref{cashinvariance}, $\tilde{v}+m\sim^* m^*_{\tilde{v}}+m$.
But this implies that
\begin{equation*}U(\tilde{v}+m)= m^*_{\tilde{v}}+m=U(\tilde{v})+m,\end{equation*}
as desired.

Next, let us show property (v). 
Assume $r\in V^\phi_0$ and $\tilde{v},\tilde{v}+r\in\tilde{V}^\phi_{0}$ with $\tilde{v},r$ being pc.
Because $v\sim^{*} m^*_{\tilde{v}}$ and $r\sim^* m^*_{r}$, we have
$$m^*_{\tilde{v}}+ m^*_{r}\leq \max_{w,s} \tilde{v}^w(s)+\max_s r(s)=\max_w\max_s (\tilde{v}^w(s)+ r(s))\in \Imp,$$
where the last equality holds by the assumed comonotonicity of $\tilde{v}$ and $r$. 
The right-hand side is in $\Imp$ because we have assumed that $\tilde{v}+r$ only takes values in $\Imp$.  
Similarly one can show that 
$m^*_{\tilde{v}}+ m^*_{r}\geq \min_w\min_s (\tilde{v}^w(s)+ r(s))\in \Imp$ and therefore we obtain as before that $m^*_{\tilde{v}}+ m^*_{r}\in \Imp.$
Hence, it follows from A7 that
$\tilde{v}+r\sim^* m^*_{\tilde{v}}+ m^*_{r}$.
Thus,
\begin{equation*}
U(\tilde{v}+r) =m^*_{\tilde{v}}+ m^*_{r}=U(\tilde{v})+U(r).
\end{equation*}
Hence, $U$ satisfies (v).

To prove property (vi), let $v\in V^\phi_{0}$ and note that $v,v$ is pc.
Thus, $U(2v)=U(v+v)=2U(v)$, by property (v) whenever $2v\in V^\phi_{0}$.
Upon iterating this argument, $U(\lambda v)=\lambda U(v)$
for all rational non-negative $\lambda$ such that $\lambda v\in V^\phi_{0}  .$
By A3, $U$ is continuous on $V_0^\phi$. 
Hence,
$U(\lambda v)=\lambda U(v)$ for all $\lambda\geq 0$ such that $\lambda v\in V^\phi_{0}$.


Property (viii) follows from Lemma~25 in \cite{MMR06}.

\smallskip
\noindent\textbf{Step 2:}\\
We prove that there exists a non-constant, non-decreasing and continuous
function $\psi:[0,1]\rightarrow [0,1]$ with $\psi(0)=0$ and $\psi(1)=1$ such that, for all $v\in V^\phi_0$,
\begin{equation}
\label{dis}
U(v)=\int_{-\infty}^0 (\psi(1-F_{v}(t))-1)\,\mathrm{d}t+\int_0^{\infty} \psi(1-F_{v}(t))\,\mathrm{d}t.
\end{equation}
Denote by $V_{0}^{[0,1]}$ all random variables in $V_0$ that only take values between $0$ and $1$. 
As $\phi(0)=0$ and $\phi(1)=1$, clearly $V_{0}^{[0,1]}$ is a subset of $V_0^\phi$. Denote $\tilde{V}_0^{[0,1]}\subset \tilde{V}_0^{\phi}$ similarly. 
The proof of Step 2 then consists of the following two parts:
\begin{itemize}
\item[(I)] First, we show that it is sufficient to prove (\ref{dis}) for $v\in V^{[0,1]}_{0}$.
\item[(II)] Then, we show that our Axiom A7, when restricted to $V_{0}$, corresponds to Axiom A7D used by \cite{Yaari87}.
Thus, we conclude that (\ref{dis}) holds.
\end{itemize}

\noindent Part (I):
It is sufficient to prove \eqref{dis} for 
$v\in V_{0}^{[0,1]} \subset  V^\phi_{0}$ because, for any $v\in   V^\phi_{0}$, there exists $a$ and $1\geq b>0$ such that $0\leq a+bv\leq 1$ (where the inequalities hold for all outcomes). 
Hence, by (vi) and (vii),
\begin{equation*}
U(v)=\frac{1}{b} U(a+bv)-\frac{a}{b}
=\frac{1}{b} \int (a+bv) \,\mathrm{d}\nu_{\psi}-\frac{a}{b}
=\int v \,\mathrm{d}\nu_{\psi}.
\end{equation*}
Therefore, it is indeed sufficient to prove (\ref{dis}) for $v\in V^{[0,1]}_0$.

\noindent Part (II):
We need the following lemma:
\begin{lemma}
\label{equ} Maintain Axioms A1-A6 (on $\tilde{V}^\phi_0$).
On the space $V^\phi_0$, Axiom A7D is then implied by A7,
i.e., for $v,u,r\in V^\phi_{0}$ with $v,r$ and $u,r$ pc and $v+r,u+r\in V^\phi_{0}$, we have
\begin{align}\label{convex}
& \mbox{for every }\alpha\in(0,1), v\succeq^*  u \Rightarrow\a v+(1-\a)r\succeq^* \a u+(1-\a)r. 
\end{align}
\end{lemma}
\textit{Proof.}
\begin{proofwithoutproofatstart}
Suppose that $v\succeq^* u$ and that the implication
\begin{align} 
\label{convex2} v\succeq^*  u \Rightarrow v+ r'\succeq^* u+r'\end{align}
holds for pc $v,r'$ and $u,r'$. 
Let $\alpha\in(0,1)$. 
Since $0\in \Imp$, we have that $v,u,r\in \Imp$ implies that also $\a v,\a u, (1-\a) r\in\Imp$.

Note that $\a v+(1-\a)r\succeq^* \a u+(1-\a)r$ would then follow directly from \eqref{convex2} if we could show that
$\a v\succeq^* \a u$ since then, 
$$\a v+(1-\a)r=\a v+ r' \succeq^* \a u+r'=\a u+(1-\a)r,$$
with $r':=(1-\a)r$.
But $\a v\succeq^* \a u$ is an immediate consequence of the fact that $U$ is a numerical representation of
$\succeq^*$ on $V^\phi_0$, and satisfies (vi) of Step 1 above.
\end{proofwithoutproofatstart}

\smallskip
Next, as in \cite{Yaari87}, we shall refer to $1-F(t)$, with $F(t)$ a CDF, as a decumulative distribution function (DDF).
We suppose in the remainder of this proof that $F$ is supported on the unit interval.
The (generalized) inverse of a DDF is a reflected (in $t=1/2$) quantile function, $F^{-1}(1-t)$.
In view of the neutrality axiom (Axiom A2),
$\succeq^*$ also induces a preference relation on the space of conditional reflected quantile functions (inverse DDFs):
given $\tilde{v}\in \tilde{V}_0^{[0,1]}$ with conditional quantile function $q_{\tilde{v}}$,
we can define its conditional reflected quantile function $\tilde{G}_{\tilde{v}}$ by
\begin{equation}
\label{revq}
\tilde{G}_{\tilde{v}}(.,t)=q_{\tilde{v}}(.,1-t),\qquad t\in[0,1].
\end{equation}
Now define
\begin{equation*}
\tilde{G}_{\tilde{v}}(\succeq^*) \tilde{G}_{\tilde{u}}\qquad\mbox{ if, and only if, }\qquad\tilde{v} \succeq^* \tilde{u}.
\end{equation*}
With this definition, we have defined a preference relation, $(\succeq^*)$, on the 
convex space
\begin{align}
\label{gamma} \tilde{\Gamma}=\{ & \tilde{G}: \,\,W\times
[0,1]\rightarrow [0,1]|\mbox{ For every fixed } t\in[0,1],\nonumber\\
&\tilde{G}(.,t)\mbox{ is }\Sigma'\mbox{-measurable. For}\mbox{ every
fixed } w\in W,\nonumber\\&\tilde{G}(w,.)\mbox{ is a decreasing and right-continuous
step function with}\nonumber\\&\tilde{G}(w,1)=0\}.
\end{align}
Indeed, every conditional reflected quantile function is in $\tilde{\Gamma}$ and for every element $G\in \tilde{\Gamma}$,
there exists a random variable $\tilde{v}\in\tilde{V}_0^{[0,1]}$ such that $G=\tilde{G}_{\tilde{v}}$.\footnote{The latter statement can be verified as follows.
Recall that if $U^{\centerdot}$ is uniformly distributed on the unit interval
and $q_X$ is the quantile function of the random variable $X$,
then $q_X(U^{\centerdot})$ has the same distribution as $X$.
Thus, if we define the random variable $\tilde{v}\in\tilde{V}_0^{[0,1]}$ through $\tilde{v}_{\tilde{G}}=\tilde{G}(w,1-U^w)$,
then the conditional reflected quantile function of $\tilde{v}$
is equal to $\tilde{G}$.
Furthermore, by neutrality, $\tilde{G}_1 (\succeq^*) \tilde{G}_2$
if, and only if, $\tilde{v}_{\tilde{G}_1}\succeq^* \tilde{v}_{\tilde{G}_2}$.}
For simplicity, we will henceforth denote the preference relations
on the spaces $\tilde{\Gamma}$ and $\tilde{V}_0 ^{[0,1]}\subset \tilde{V}_0 ^\phi$
both by $\succeq^*$, too.
We define $\Gamma$ as the subspace of all elements in $\tilde{\Gamma}$ that carry no
ambiguity, i.e.,
\begin{align*}
\Gamma=\{ & G\in\tilde{\Gamma}| \,\,\mbox{for all }w_1,w_2\in W:\,\,\,
G(w_1,\cdot)=G(w_2,\cdot)\}.
\end{align*}

Lemma \ref{equ} implies that, on the space of non-negative random
variables in $V_0$ bounded by one, Axioms A1-A4 and A7D hold.
Similar to \cite{Yaari87}, by neutrality (A2), this induces
a preference relation on the space of
conditional reflected quantile functions, $\Gamma$,
simply denoted by $\succeq^*$, that satisfies
weak and non-degenerate order, continuity,
certainty first-order stochastic dominance
and the independence axiom.\footnote{The independence axiom asserts that if, for DDFs $G_1,G_2,G_3\in \Gamma$, $G_1\succeq^* G_2$,
then, for every $\a\in(0,1)$, $\a G_1+(1-\a)G_3\succeq^* \a G_2+(1-\a)G_3$.}
Therefore, by the mixture space theorem (\citealp{HersteinMilnor53}),
there exists a non-constant, non-decreasing and continuous function $\psi:[0,1]\rightarrow [0,1]$
such that the corresponding numerical representation is given by
\begin{equation*}U(v)=-\int_0^1 \psi(t)G_{v}(.,\,\mathrm{d}t)
=\int_0^1 \psi(G^{-1}_{v}(.,t))\,\mathrm{d}t
=\int_0^1 \psi(1-F_{v}(.,t))\,\mathrm{d}t,
\end{equation*}
where $G_{v}$ is defined by (\ref{revq}).
Hence, (\ref{dis}) holds.
Finally, it is straightforward to verify that $U(m)=m$ for all $m\in[0,1]$ (see Step 1 above)
implies that
we must have $\psi(0)=0$ and $\psi(1)=1$.

\smallskip
\noindent\textbf{Steps 3+4:}\\
Recall Step 1.
By construction, $U(0)=0$.
As $U$ satisfies (i)-(viii), 
$U$ may be identified with a concave and normalized
niveloid\footnote{The mapping $U$ from $V^{\prime,\phi}_0$ to $\Imp$ is a concave and normalized niveloid
if it is concave, Lipschitz continuous with respect to the $||.||_\infty$-norm, and satisfies $U(m)=m$ for all $m\in\Imp$.} on the space
of bounded, $\Sigma'$-measurable functions on $W$; see Lemma~25 in \cite{MMR06}.
Duality results in convex analysis for niveloids
(see Lemma 26 in \citealp{MMR06}) then yield that, for all $v'\in V^{\prime,\phi}_0$,
\begin{equation}\label{xprime}
U(v')=\min_{Q\in \Delta(W,\Sigma')}  \left\{\E_Q\left[v'\right]+c_{\mathrm{min}}(Q)\right\},
\end{equation}
with $c_{\mathrm{min}}$ defined by
\begin{align}
\label{defcmin}
c_{\mathrm{min}}(Q)=\sup_{v'\in V^{\prime,\phi}_0}\{U(v')-\E_Q\left[v'\right]\}\geq U(0)=0,
\end{align}
and being the unique minimal function satisfying (\ref{xprime}).
As $U(m)=m$ for all $m\in\Imp$, there exists a $Q$ such that $c(Q)<\infty$.
Now we have
\begin{equation*}0=U(0)=\min_{Q\in \Delta(W,\Sigma')} c(Q).\end{equation*}
In particular, $c$ is grounded, convex and lower-semicontinuous.

For $\tilde{v}\in\tilde{V}^\phi_0$, define $m^*_{\tilde{v}^{w}}\in \Imp$
as the corresponding \emph{certainty equivalent} of $\tilde{v}$
in the state of the world $w$, i.e.,
\begin{equation*}m^*_{\tilde{v}^{w}}=U(\tilde{v}^{w})=\int  \tilde{v}^w \,\mathrm{d}\nu_{\psi};\end{equation*}
see Lemma~\ref{certequ}.
Set $\bar{v}^w=m^*_{\tilde{v}^{w}}$.
Clearly, $\bar{v}$ is independent of $s$.
Furthermore, by the Theorem of Tornelli, $\bar{v}$ is
$\Sigma'$-measurable. 
In particular, $\bar{v}$ is in $V^{\prime,\phi}_0$.
Observe that, for every $w\in W$,
$U(\bar{v}^{w})=U(m^*_{\tilde{v}^{w}})=U(U(\tilde{v}^{w}))=U(\tilde{v}^{w})$,
where we have used in the last equality that, for all $m\in\Imp$, $U(m)=m$ and by the monotonicity property~(iv) of Step~1,
\begin{align*}
\Imp\ni\min_s \tilde{v}^w(s) &=U(\min_s \tilde{v}^w(s))\leq U( \tilde{v}^w)\\ 
&\leq U(\max_s \tilde{v}^w(s)) =\max_s \tilde{v}^w(s)\in \Imp,
\end{align*}
so that $U( \tilde{v}^w)\in\Imp$.
Hence, property (iv) of Step~1 implies that $U(\bar{v})=U(\tilde{v})$.
This entails that, for all $\tilde{v}\in \tilde{V}^\phi_{0}$,
\begin{align*}
U(\tilde{v})=&\ U(\bar{v})
=
\min_{Q\in \Delta(W,\Sigma')}\left\{\E_Q[\bar{v}]+c(Q)\right\}\\
=&\min_{Q\in \Delta(W,\Sigma')}\left\{\int U(m^*_{\tilde{v}^{w}})Q(\mathrm{d}w)+c(Q)\right\}\\
=&\min_{Q\in \Delta(W,\Sigma')}\left\{\int U(\tilde{v}^{w})Q(\mathrm{d}w)+c(Q)\right\}\\
=&\min_{Q\in \Delta(W,\Sigma')}\left\{\int\left(\int_{-\infty}^0
(\psi(1-F_{\tilde{v}^w}(t))-1)\,\mathrm{d}t\nonumber\right.\right.\\
&\qquad\qquad\qquad\left.\left.+\int_0^{\infty} \psi(1-F_{\tilde{v}^w}(t))\,\mathrm{d}t\right)Q(\mathrm{d}w)+c(Q)\right\}\\
=&\min_{Q \in \Delta(W,\Sigma')}\left\{\E_Q\left[\int \tilde{v}^{\centerdot} \,\mathrm{d}\nu_{\psi}\right]+c(Q)\right\},
\end{align*}
where we have used (\ref{xprime}) in the second and (\ref{dis})
in the fifth equalities.
This proves the ``only if'' part of Theorem \ref{1}($\alpha$).

The proof of Theorem \ref{1}($\beta$) now 
follows by defining,
for all $\tilde{v},\tilde{u}\in\tilde{V}^\phi$,
\begin{align*}
\tilde{v}\succeq^* \tilde{u}\Leftrightarrow
U(\tilde{v}):&=\min_{Q \in \Delta(W,\Sigma')}\left\{\E_Q\left[\int \tilde{v}^{\centerdot}\,\mathrm{d}\nu_{\psi}\right]+c(Q)\right\}\\
&\geq \min_{Q \in \Delta(W,\Sigma')}\left\{\E_Q\left[\int \tilde{u}^{\centerdot}\,\mathrm{d}\nu_{\psi}\right]+c(Q)\right\}=U(\tilde{u}).
\end{align*}
The extension is unique, as it may be seen as before that it follows from our axioms that the functional $U$ defined above is concave on $\tilde{V}^\phi$. 
Therefore, $U$ is continuous on the interior of its domain, which gives a unique extension from the dense subspace $\tilde{V}^\phi_0$ to $\tilde{V}^\phi$.
\end{proofwithoutproofatstart}

\smallskip
\textit{Proof of Theorem~\ref{th:GS}.}
\begin{proofwithoutproofatstart}
For the proof of Theorem~\ref{th:GS} we need the following lemma:
\begin{lemma}
\label{scaling}
A7$^{0}$ with ``$\oplus$'' replaced by ``+'' for random variables only taking values in $\Imp$ (i.e., with $\tilde{V}_0,V_0$ replaced by $\tilde{V}^\phi_0,V^\phi_0$) implies that, for $\tilde{v},\tilde{u}\in\tilde{V}^\phi_0$,
$\tilde{v}\succeq^* \tilde{u}$ if, and only if,
$\lambda \tilde{v}\succeq^* \lambda\tilde{u}$ for every $\lambda\geq 0$ such that $\lambda \tilde{v},\lambda \tilde{u}\in\tilde{V}^\phi_0$.
\end{lemma}
\textit{Proof.}
\begin{proofwithoutproofatstart}
The proof of the ``if'' part is straightforward.
Let us prove the ``only if'' part.
So, suppose that $\tilde{v}\succeq^* \tilde{u}$.
If $\lambda \in[0,1]$, then $\lambda \tilde{v}\succeq^* \lambda\tilde{u}$ follows
directly from Axiom A7$^{0}$ with $\a=\lambda$ and $r=0$.
If $\lambda>1$, then let us suppose that $\lambda \tilde{u}\succ^* \lambda\tilde{v}$ would hold.
Defining $\a=\frac{1}{\lambda}\in(0,1)$ yields, by A7$^{0}$,
\begin{equation*}\tilde{u}=\a \lambda\tilde{u}+(1-\a)0\succ^* \a \lambda\tilde{v}+(1-\a)0=\tilde{v},\end{equation*}
which is a contradiction.
Hence, indeed $\lambda \tilde{v}\succeq^* \lambda\tilde{u}$ for every $\lambda\geq 0$.
\end{proofwithoutproofatstart}
\smallskip
Thus, Axiom A7$^{0}$ implies that the preference relation is scale invariant on $\tilde{V}^\phi_0$,
whereas Axiom A7 only implies scale invariance on $V^\phi_0$ (formally, via Lemma~\ref{equ}).
The next proposition shows explicitly that Axiom A7$^{0}$ is stronger than Axiom A7:

\begin{proposition}\label{prop:A70}
Axiom A7$^{0}$ implies Axiom A7, with ``$\oplus$'' replaced by ``+'' for random variables only taking values in $\Imp$ (i.e., with $\tilde{V}_0$, $V_0$ replaced by $\tilde{V}_0^\phi$, $V_0^\phi$).
\end{proposition}
\textit{Proof.}
\begin{proofwithoutproofatstart}
First, in view of Lemma~\ref{scaling}, we can extend $\succeq^*$ consistently to $\tilde{V}_0$ by defining, for \emph{any} $\tilde{v},\tilde{u}\in \tilde{V}_0$, that
$\tilde{v}\succeq^*\tilde{u}$ if
$$\frac{\min(-\inf \phi,1)}{2\max(|\tilde{v}|\vee |\tilde{u}|)}\tilde{v}\succeq^* \frac{\min(-\inf \phi,1)}{2\max(|\tilde{v}|\vee |\tilde{u}|)}\tilde{u},$$
where in case both random variables are degenerate at zero the scaling factor should be set equal to $1$.
Note that as $0,1\in \Imp$ and $\Imp$ is open, we have $\inf\phi<0$, and therefore
$\inf\phi< \frac{\min(-\inf \phi,1)}{2\max(|\tilde{v}|\vee |\tilde{u}|)}\tilde{v}< 1$
and a similar inequality holds for the scaled $\tilde{u}$, such that the scaled random variables only take values in $\Imp$. 
Clearly, the (scaled) extension of $\succeq^*$ to $\tilde{V}_0$ is again scale invariant in the sense that, \textit{mutatis mutandis}, Lemma~\ref{scaling} holds on the entire space $\tilde{V}_0$. 
Thus, Axiom A7$^0$ also holds on the entire space $\tilde{V}_0$.

Next, let us show the proposition. 
Suppose that $\tilde{v}\succeq^* \tilde{u}$ and that $\tilde{v},\tilde{u}$ and $r$ are pc.
Let $\a\in(0,1).$
Then, by Lemma~\ref{scaling} and the paragraph above, under Axiom A7$^{0}$, $\frac{1}{\a}\tilde{v}\succeq^* \frac{1}{\a}\tilde{u}$.
Next, let $\bar{r}=\frac{r}{1-\a}$.
Then, we obtain from Axiom A7$^{0}$ that
\begin{equation*}\tilde{v}+r=\a \left(\frac{1}{\a}\tilde{v}\right)+(1-\a)\bar{r}\succeq^* \a \left(\frac{1}{\a}\tilde{u}\right)+(1-\a)\bar{r}=
\tilde{u}+r.\end{equation*}
Hence, A7 is indeed satisfied.
\end{proofwithoutproofatstart}
Because Axioms A1-A6 and A7$^0$ imply A1-A7, we may, as before, obtain a numerical representation $U$ for $\succeq^*$. 
Noticing that, by Lemma~\ref{scaling}, $U$ must be positively homogeneous, we can as in the proof of Proposition~\ref{prop:A70} extend $U$ to the entire space $\tilde{V}_0$, and it follows from classical results in convex analysis that the minimal function $c_{\mathrm{min}}$ in \eqref{xprime} can then be chosen to only take the values zero or infinity. 

The unique extension of $\succeq^*$ to $\tilde{V}$ can be seen as in the proof of Theorem~\ref{1}.  
The proof of Theorem~\ref{th:GS}(b) follows from \cite{FollmerSchied16}, Chapter 4.
This completes the proof of Theorem~\ref{th:GS}.
\end{proofwithoutproofatstart}

\textit{Proof of Proposition~\ref{prop:SSD}.}
\begin{proofwithoutproofatstart}
Recall that $v \succeq^* w$ if, and only if, $\phi^{-1}(v) \succeq \phi^{-1}(w)$, with $\phi$ from Proposition~\ref{prop:app}. 
Furthermore, recall the equivalence asserted by Lemma~\ref{starpreference}. 
Now one may verify that the corresponding Axioms A1-A3 are already sufficient
to guarantee the existence of a numerical representation of $\succeq^*$. 
As in the proof of Theorem~\ref{th:RDU}, denote it by $U$.
By~\eqref{Urep}, $U$ may be seen to be concave and upper-semicontinuous on the space $V^{\prime,\phi}_0$. 
Therefore, Theorem~4.1 in \cite{Dana05} implies that $\succeq^*$ respecting SSD is equivalent to $U$ being law invariant under $P'$. 
Using the definition of $\succeq^*$ then shows the first part of the proposition.  
Noting that, if $\phi$ is concave, $v'$ dominates $u'$ in SSD if, and only if, $\phi(u')$ dominates $\phi(v')$ in SSD finishes the proof of Proposition~\ref{prop:SSD}.
\end{proofwithoutproofatstart}

\textit{Proof of Proposition~\ref{prop:caa}.}
\begin{proofwithoutproofatstart}
If $\succeq^*$ is more ambiguity averse than $\succeq$,
then $\succeq$ and $\succeq^*$ agree on $V_0$, and therefore we may choose a positive affine transformation such that $\phi^* = \phi$, and furthermore $\psi^* = \psi$. 
This implies that 
\begin{equation*}
c^*_{\mathrm{min}}(Q)=\sup_{v'\in V'_0}\left\{m^*_{\phi(v')}-\E_Q[\phi(v')]\right\}
\geq \sup_{v'\in V'_0}\left\{m_{\phi(v')}-\E_Q[\phi(v')]\right\}=c_{\mathrm{min}}(Q).
\end{equation*}
This proves the ``only if'' part.
To prove the ``if'' part, suppose that $c^*\geq c$, $\psi^*=\psi$ and $\phi^* = \phi$.
Then $\tilde{v}\succeq v$, with $\tilde{v}\in\tilde{V}_{0}$ and $v\in V_{0}$, entails that
\begin{align*}
m_{v}^*=m_{v}
&\leq \phi^{-1}\Big(\min_{Q\in\Delta(W,\Sigma')}\left\{\E_Q\left[\int \phi(\tilde{v}^{\centerdot}) \,\mathrm{d}\nu_{\psi}\right]+c_{\mathrm{min}}(Q)\right\}\Big)\\
&\leq \phi^{*,-1}\Big(\min_{Q\in\Delta(W,\Sigma')}\left\{\E_Q\left[\int \phi^{*}(\tilde{v}^{\centerdot}) \,\mathrm{d}\nu_{\psi^*}\right]+c^*_{\mathrm{min}}(Q)\right\}\Big)=m_{\tilde{v}}^*.
\end{align*}
Hence, $\tilde{v}\succeq^* v$.
\end{proofwithoutproofatstart}

\textit{Proof of Proposition~\ref{prop:aa}.}
\begin{proofwithoutproofatstart}
Assume Axioms A1-A8.
By Theorem~\ref{th:RDU}, there exist functions $\phi$, $\psi$ and $c$ such that \eqref{eq:rduamb} holds.
Set $P'=\argmin_Q c_{\mathrm{min}}(Q)$.
Because $c_{\mathrm{min}}$ is grounded, $c_{\mathrm{min}}(P')=0$.
Denote by $\succeq^{\mathrm{AN}}$ the ambiguity neutral agent with measure $P'$, utility function $\phi$ and probability weighting function $\psi$.
Suppose that $\tilde{v}\succeq v$, $\tilde{v}\in \tilde{V}_0$ and
$v\in V_0$.
Then,
\begin{align*}m_{v}&\leq  \phi^{-1}\Big(\min_{Q\in \Delta(W,\Sigma')}\left\{ \mathbb{E}_{Q}\left[
\int \phi(\tilde{v}^{\centerdot}) \,\mathrm{d}\nu_{\psi}\right]+c_{\mathrm{min}}(Q)\right\}\Big)\\&\leq   \phi^{-1}\Big(\mathbb{E}_{P'}\left[\int
\phi(\tilde{v}^{\centerdot}) \,\mathrm{d}\nu_{\psi}\right]+c_{\mathrm{min}}(P')\Big)
=
\phi^{-1}\Big(\mathbb{E}_{P'}\left[ \int \phi(\tilde{v}^{\centerdot}) \,\mathrm{d}\nu_{\psi}\right]\Big).
\end{align*}
\end{proofwithoutproofatstart}

\end{appendix}

%



{\small

\begin{thebibliography}{99}

\baselineskip 0.37 cm


\bibitem[\protect\citeauthoryear{Allais}{1953}]{Allais53}
\textsc{Allais, M.} (1953).
Le comportement de l'homme rationnel
devant le risque: critique des postulats et axiomes de l'\'ecole Am\'ericaine.
\textit{Econometrica} \textbf{21}, 503--546.


\bibitem[\protect\citeauthoryear{Anscombe and Aumann}{1963}]{AnscombeAumann63}
\textsc{Anscombe, F. J., and Aumann, R. J.} (1963).
A definition of subjective probability.
\textit{The Annals of  Mathematical Statistics} \textbf{34}, 199--205.

\bibitem[\protect\citeauthoryear{Artzner et al.}{1999}]{ADEH99}
\textsc{Artzner, P., Delbaen, F., Eber, J., and Heath, D.} (1999).
Coherent measures of risk.
\textit{Mathematical Finance} \textbf{9}, 203--228.

\bibitem[\protect\citeauthoryear{Basak and Shapiro}{2001}]{BS01}
\textsc{Basak, S., and Shapiro, A.} (2001).
Value-at-risk-based risk management: Optimal policies and asset prices.
\textit{The Review of Financial Studies} \textbf{14}, 371--405.




\bibitem[\protect\citeauthoryear{Ben-Tal}{1985}]{Ben-Tal85}
\textsc{Ben-Tal, A.} (1985).
The entropic penalty approach to stochastic programming.
\textit{Mathematics of Operations Research} \textbf{10}, 224--240.

\bibitem[\protect\citeauthoryear{Berger and Eeckhoudt}{2021}]{BE21}
\textsc{Berger, L., and Eeckhoudt, L. R.} (2021).
Risk, ambiguity, and the value of diversification.
\textit{Management Science} \textbf{67}, 1639--1647.

\bibitem[\protect\citeauthoryear{Bernoulli}{1738}]{Bernoulli38}
\textsc{Bernoulli, D.} (1738).
Specimen theoriae novae de mensura sortis.
In: Commentarii Academiae Scientiarum Imperialis Petropolitannae (1738).
Translated from Latin into English by Sommer, L. (1954).
Exposition of a new theory on the measurement of risk.
\textit{Econometrica} \textbf{22}, 23--36.

\bibitem[\protect\citeauthoryear{Cai, Li and Mao}{2023}]{CLM22}
\textsc{Cai, J., Li, J. Y.-M. and Mao, T.} (2024).
Distributionally robust optimization under distorted expectations,
\textit{Operations Research}, in press.

\bibitem[\protect\citeauthoryear{Carr, Geman and Madan}{2001}]{PCh01}
\textsc{Carr, P., Geman, H., and Madan, D. B.} (2001).
Pricing and hedging in incomplete markets.
\textit{Journal of Financial Economics} \textbf{62}, 131--167.


\bibitem[\protect\citeauthoryear{Cherny}{2006}]{C06}
\textsc{Cherny, A.} (2006).
Weighted V@R and its properties.
\textit{Finance and Stochastics} \textbf{10}, 367--393.

\bibitem[\protect\citeauthoryear{Cherny and Madan}{2009}]{CM09}
\textsc{Cherny, A., and Madan, D. B.} (2009).
New measures for performance evaluation.
\textit{Review of Financial Studies} \textbf{22}, 2571--2606.


\bibitem[\protect\citeauthoryear{Chew and Karni}{1994}]{ChewKarni94}
\textsc{Chew, S. H., and Karni, E.} (1994).
Choquet expected utility with a finite state space: Commutativity and act-independence.
\textit{Journal of Economic Theory} \textbf{62}, 469--479.




\bibitem[\protect\citeauthoryear{Csisz\'ar}{1975}]{Csiszar75}
\textsc{Csisz\'ar, I.} (1975).
$I$-divergence geometry of probability distributions and minimization problems.
\textit{The Annals of Probability} \textbf{3}, 146--158.

\bibitem[\protect\citeauthoryear{Dana}{2005}]{Dana05}
\textsc{Dana, R.-A.} (2005).
A representation result for concave Schur concave functions.
\textit{Mathematical Finance} \textbf{15}, 613--634.

\bibitem[\protect\citeauthoryear{de Finetti}{1931}]{deFinetti31}
\textsc{de Finetti, B.} (1931).
Sul significato soggettivo della probabilit\`a.
\textit{Fundamenta Mathematicae} \textbf{17}, 298--329.

\bibitem[\protect\citeauthoryear{Dean and Ortoleva}{2017}]{DeanOrtoleva17}
\textsc{Dean, M., and Ortoleva, P.} (2017).
Allais, Ellsberg, and preferences for hedging.
\textit{Theoretical Economics} \textbf{12}, 377--424.


\bibitem[\protect\citeauthoryear{Denneberg}{1994}]{D94}
\textsc{Denneberg, D.} (1994).
\textit{Non-Additive Measure and Integral.}
Springer, Dordrecht.

\bibitem[\protect\citeauthoryear{Duffie and Pan}{1997}]{DP97}
\textsc{Duffie, D., and Pan, J.} (1997).
An overview of value at risk.
\textit{The Journal of Derivatives} \textbf{4}, 7--49.

\bibitem[\protect\citeauthoryear{Eeckhoudt, Laeven and Schlesinger}{2020}]{ELS20}
\textsc{Eeckhoudt, L. R., Laeven, R. J. A., and Schlesinger, H.} (2020).
Risk apportionment: {T}he dual story.
\textit{Journal of Economic Theory}, 104971.

\bibitem[\protect\citeauthoryear{Eeckhoudt and Laeven}{2022}]{EL21}
\textsc{Eeckhoudt, L. R., and Laeven, R. J. A.} (2022).
Dual moments and risk attitudes.
\textit{Operations Research} \textbf{70}, 1330--1341.

\bibitem[\protect\citeauthoryear{Ellsberg}{1961}]{Ellsberg61}
\textsc{Ellsberg, D.} (1961).
Risk, ambiguity and the Savage axioms.
\textit{Quarterly Journal of Economics} \textbf{75}, 643--669.

\bibitem[\protect\citeauthoryear{Embrechts, Puccetti, and Rüschendorf}{2013}]{Embrechts13}
\textsc{Embrechts, P., Puccetti, G., and Rüschendorf, L.} (2013).
Model uncertainty and VaR aggregation.
\textit{Journal of Banking and Finance} \textbf{37}, 2750--2764.

\bibitem[\protect\citeauthoryear{Epstein}{1999}]{Epstein99}
\textsc{Epstein, L. G.} (1999).
A definition of uncertainty aversion.
\textit{Review of Economic Studies} \textbf{66}, 579--608.


\bibitem[\protect\citeauthoryear{F\"ollmer and Schied}{2002}]{FS02}
\textsc{F\"ollmer, H., and Schied, A.} (2002).
Convex measures of risk and trading constraints.
\textit{Finance and Stochastics} \textbf{6}, 429--447.

\bibitem[\protect\citeauthoryear{F\"ollmer and Schied}{2016}]{FollmerSchied16}
\textsc{F\"ollmer, H., and Schied, A.} (2016).
\textit{Stochastic Finance}.
2nd ed., De Gruyter, Berlin (4th ed. 2016).

\bibitem[\protect\citeauthoryear{Frittelli and Rosazza Gianin}{2002}]{FRG02}
\textsc{Frittelli, M., and Rosazza Gianin, E.} (2002).
Putting order in risk measures.
\textit{Journal of Banking \& Finance} \textbf{26}, 1473--1486.




\bibitem[\protect\citeauthoryear{Ghaoui, Oks and Oustry}{2003}]{GOO03}
\textsc{Ghaoui, L. E., Oks, M., and Oustry, F.} (2003). 
Worst-case value-at-risk and robust portfolio optimization: A conic programming approach. 
\textit{Operations Research} \textbf{51}, 543--556.

\bibitem[\protect\citeauthoryear{Ghirardato and Marinacci}{2001}]{GhirardatoMarinacci01}
\textsc{Ghirardato, P., and Marinacci, M.} (2001).
Risk, ambiguity, and the separation of utility and beliefs.
\textit{Mathematics of Operations Research} \textbf{26}, 864--890.

\bibitem[\protect\citeauthoryear{Ghirardato and Marinacci}{2002}]{GhirardatoMarinacci02}
\textsc{Ghirardato, P., and Marinacci, M.} (2002).
Ambiguity made precise: A comparative foundation.
\textit{Journal of Economic Theory} \textbf{102}, 251--289.

\bibitem[\protect\citeauthoryear{Ghirardato, Maccheroni, Marinacci and Siniscalchi}{2003}]{GhirardatoMaccheroniMarinacciSiniscalchi03}
\textsc{Ghirardato, P., Maccheroni, F., Marinacci, M., and Siniscalchi, M.} (2003).
A subjective spin on roulette wheels.
\textit{Econometrica} \textbf{71}, 1897--1908.

\bibitem[\protect\citeauthoryear{Ghirardato, Maccheroni and Marinacci}{2004}]{GMM04}
\textsc{Ghirardato, P., Maccheroni, F., and Marinacci, M.} (2004). 
Differentiating ambiguity and ambiguity attitude. 
\textit{Journal of Economic Theory} \textbf{118}, 133--173.

\bibitem[\protect\citeauthoryear{Gilboa and Schmeidler}{1989}]{GilboaSchmeidler89}
\textsc{Gilboa, I., and Schmeidler, D.} (1989).
Maxmin expected utility with non-unique prior.
\textit{Journal of Mathematical Economics} \textbf{18}, 141--153.

\bibitem[\protect\citeauthoryear{Gonzalez and Wu}{1999}]{GonzalezWu99}
\textsc{Gonzalez, R., and Wu, G.} (1999).
On the shape of the probability weighting function.
\textit{Cognitive Psychology} \textbf{38}, 129--166.


\bibitem[\protect\citeauthoryear{Gul}{1992}]{Gul92}
\textsc{Gul, F.} (1992).
Savage's theorem with a finite number of states.
\textit{Journal of Economic Theory} \textbf{57}, 99--110.

\bibitem[\protect\citeauthoryear{Gul and Pesendorfer}{2015}]{GulPesendorfer15}
\textsc{Gul, F., and Pesendorfer, W.} (2015). 
Hurwicz expected utility and subjective sources. 
\textit{Journal of Economic Theory} \textbf{159}, 465--488.


\bibitem[\protect\citeauthoryear{Hansen and Sargent}{1995}]{HS95}
\textsc{Hansen, L. P., and Sargent. T. J.} (1995).
Discounted linear exponential quadratic Gaussian control.
\textit{IEEE Transactions on Automatic Control} \textbf{40}, 968-971.

\bibitem[\protect\citeauthoryear{Hansen and Sargent}{2000}]{HansenSargent00}
\textsc{Hansen, L. P., and Sargent, T. J.} (2000).
Wanting robustness in macroeconomics.
Mimeo, University of Chicago and Stanford University.

\bibitem[\protect\citeauthoryear{Hansen and Sargent}{2001}]{HS01}
\textsc{Hansen, L. P., and Sargent, T. J.} (2001).
Robust control and model uncertainty.
\textit{American Economic Review} \textbf{91}, 60--66.

\bibitem[\protect\citeauthoryear{Hansen and Sargent}{2007}]{HS07}
\textsc{Hansen, L. P., and Sargent, T. J.} (2007).
\textit{Robustness.}
Princeton University Press, Princeton.

\bibitem[\protect\citeauthoryear{Hansen}{2014}]{H14}
\textsc{Hansen, L. P.} (2014).
Nobel lecture: {U}ncertainty outside and inside economic models.
\textit{Journal of Political Economy} \textbf{122}, 945--987.


\bibitem[\protect\citeauthoryear{Harrison and Swarthout}{2016}]{HarrisonSwarthout16}
\textsc{Harrison, G. W., and Swarthout, J. T.} (2016).
Cumulative prospect theory in the laboratory: {A} reconsideration.
Mimeo, CEAR.

\bibitem[\protect\citeauthoryear{He, Jin and Zhou}{2015}]{HX15}
\textsc{He, X. D., Jin, H., and Zhou, X. Y.} (2015).
Dynamic portfolio choice when risk is measured by weighted VaR.
\textit{Mathematics of Operations Research} \textbf{40}, 773-796.


\bibitem[\protect\citeauthoryear{Herstein and Milnor}{1953}]{HersteinMilnor53}
\textsc{Herstein, I. N. and J. Milnor} (1953).
An axiomatic approach to measurable utility.
\textit{Econometrica} \textbf{21}, 291--297.

\bibitem[\protect\citeauthoryear{Huber}{1981}]{Huber81}
\textsc{Huber, P. J.} (1981).
\textit{Robust Statistics}.
Wiley, New York.


\bibitem[\protect\citeauthoryear{Jorion}{1997}]{J97}
\textsc{Jorion, P.} (1997).
\textit{Value at Risk}.
McGraw-Hill, New York.


\bibitem[\protect\citeauthoryear{Keynes}{1921}]{Keynes21}
\textsc{Keynes, J. M.} (1921).
\textit{A Treatise on Probability}.
Macmillan and Company Ltd., London.


\bibitem[\protect\citeauthoryear{Knight}{1921}]{Knight21}
\textsc{Knight, F. H.} (1921).
\textit{Risk, Uncertainty, and Profit}.
Houghton Mifflin, Boston.


\bibitem[\protect\citeauthoryear{Kou, Peng and Heyde}{2013}]{KPH13}
\textsc{Kou, S., Peng, X., and Heyde, C. C.} (2013).
External risk measures and Basel accords.
\textit{Mathematics of Operations Research} \textbf{38}, 393-417.

\bibitem[\protect\citeauthoryear{Kusuoka}{2001}]{K01}
\textsc{Kusuoka, S.} (2001).
On law invariant coherent risk measures.
\textit{Advances in Mathematical Economics}
\textbf{3}, 83-95.

\bibitem[\protect\citeauthoryear{Laeven and Stadje}{2013}]{LaevenStadje13}
\textsc{Laeven, R. J. A., and Stadje, M.} (2013).
Entropy coherent and entropy convex measures of risk.
\textit{Mathematics of Operations Research} \textbf{38}, 265--293.

\bibitem[\protect\citeauthoryear{Laeven and Stadje}{2014}]{LaevenStadje14}
\textsc{Laeven, R. J. A., and Stadje, M.} (2014).
Robust portfolio choice and indifference valuation.
\textit{Mathematics of Operations Research} \textbf{39}, 1109--1141.

\bibitem[\protect\citeauthoryear{Lesnevski, Nelson and Staum}{2007}]{LNS07}
\textsc{Lesnevski, V., Nelson, B., and Staum, J.} (2007).
Simulation of coherent risk measures based on generalized scenarios.
\textit{Management Science} \textbf{53}, 1756--1769.

\bibitem[\protect\citeauthoryear{Li}{2018}]{L22} 
\textsc{Li, J. Y. M.}(2018)
Closed-form solutions for worst-case law invariant risk measures with application to robust portfolio optimization. 
\textit{Operations Research} \textbf{66}, 1533--1541.

\bibitem[\protect\citeauthoryear{Maccheroni}{2004}]{Maccheroni04}
\textsc{Maccheroni, F.} (2004).
Yaari's dual theory without the completeness axiom.
\textit{Journal of Economic Theory} \textbf{23}, 701--714.

\bibitem[\protect\citeauthoryear{Maccheroni, Marinacci, Rustichini and Taboga}{2004}]{MMRT04}
\textsc{Maccheroni, F., Marinacci, M., Rustichini, A., and Taboga, M.} (2004).
Portfolio selection with monotone mean-variance preferences.
ICER Working Papers - Applied Mathematics Series 27-2004,
International Centre for Economic Research.

\bibitem[\protect\citeauthoryear{Maccheroni, Marinacci and Rustichini}{2006}]{MMR06}
\textsc{Maccheroni, F., Marinacci, M., and Rustichini, A.} (2006).
Ambiguity aversion, robustness, and the variational representation of preferences.
\textit{Econometrica} \textbf{74}, 1447--1498.

\bibitem[\protect\citeauthoryear{Machina}{1987}]{Machina87}
\textsc{Machina, M. J.} (1987).
Choice under uncertainty: Problems solved and unsolved.
\textit{Journal of Economic Perspectives} \textbf{1}, 121--154.

\bibitem[\protect\citeauthoryear{Machina and Schmeidler}{1992}]{MachinaSchmeidler92}
\textsc{Machina, M. J., and Schmeidler, D.} (1992).
A more robust definition of subjective probability.
\textit{Econometrica} \textbf{60}, 745--780.

\bibitem[\protect\citeauthoryear{Machina}{2009}]{Machina09}
\textsc{Machina, M. J.} (2009).
Risk, ambiguity, and the rank-dependence axioms.
\textit{American Economic Review} \textbf{99}, 385--392.



\bibitem[\protect\citeauthoryear{Marinacci}{2002}]{Marinacci02}
\textsc{Marinacci, M.} (2002).
Probabilistic sophistication and multiple priors.
\textit{Econometrica} \textbf{70}, 755--764.


\bibitem[\protect\citeauthoryear{Nakamura}{1990}]{Nakamura90}
\textsc{Nakamura, Y.} (1990).
Subjective expected utility with non-additive probabilities on finite state spaces.
\textit{Journal of Economic Theory} \textbf{51}, 346--366.

\bibitem[\protect\citeauthoryear{Pesenti, Wang and Wang}{2022}]{PWW22}
\textsc{Pesenti, S. M., Wang, Q., and Wang, R.} (2022).
Optimizing distortion riskmetrics with distributional uncertainty.
ArXiv.

\bibitem[\protect\citeauthoryear{Prelec}{1998}]{Prelec98}
\textsc{Prelec, D.} (1998).
The probability weighting function.
\textit{Econometrica} \textbf{66}, 497--527.

\bibitem[\protect\citeauthoryear{Quiggin}{1982}]{Quiggin82}
\textsc{Quiggin, J.} (1982).
A theory of anticipated utility.
\textit{Journal of Economic Behaviour and Organization} \textbf{3}, 323--343.


\bibitem[\protect\citeauthoryear{Ramsey}{1931}]{Ramsey31}
\textsc{Ramsey, F. P.} (1931).
\textit{Truth and probability}.
In: \textit{The Foundations of Mathematics and Other Logical Essays}.
Routledge and Kegan Paul, London.

\bibitem[\protect\citeauthoryear{Riedel}{2009}]{Riedel09}
\textsc{Riedel, F.} (2009).
Optimal stopping with multiple priors.
\textit{Econometrica} \textbf{77}, 857--908.


\bibitem[\protect\citeauthoryear{Rothschild and Stiglitz}{1970}]{RothschildStiglitz70}
\textsc{Rothschild, M., and Stiglitz, J. E.} (1970).
Increasing risk: I. A definition.
\textit{Journal of Economic Theory} \textbf{2}, 225--243.

\bibitem[\protect\citeauthoryear{Ruszczy\'nsky and Vanderbei}{2003}]{RV03}
\textsc{Ruszczy\'nski, A., and Vanderbei, R. J.} (2003).
Frontiers of stochastically nondominated portfolios.
\textit{Econometrica} \textbf{71}, 1287--1297.

\bibitem[\protect\citeauthoryear{Ruszczy\'nski and Shapiro}{2006}]{RS06a}
\textsc{Ruszczy\'nski, A., and Shapiro, A.} (2006).
Optimization of convex risk functions.
\textit{Mathematics of Operations Research} \textbf{31}, 433--452.

\bibitem[\protect\citeauthoryear{Savage}{1954}]{Savage54}
\textsc{Savage, L. J.} (1954).
\textit{The Foundations of Statistics}.
Wiley, New York (2nd ed. 1972, Dover, New York).

\bibitem[\protect\citeauthoryear{Schmeidler}{1986}]{Schmeidler86}
\textsc{Schmeidler, D.} (1986).
Integral representation without additivity.
\textit{Proceedings of the American Mathematical Society} \textbf{97}, 253--261.

\bibitem[\protect\citeauthoryear{Schmeidler}{1989}]{Schmeidler89}
\textsc{Schmeidler, D.} (1989).
Subjective probability and expected utility without additivity.
\textit{Econometrica} \textbf{57}, 571--587.

\bibitem[\protect\citeauthoryear{Strzalecki}{2011a}]{Strzalecki11a}
\textsc{Strzalecki, T.} (2011a).
Axiomatic foundations of multiplier preferences.
\textit{Econometrica} \textbf{79}, 47--73.

\bibitem[\protect\citeauthoryear{Strzalecki}{2011b}]{Strzalecki11b}
\textsc{Strzalecki, T.} (2011b).
Probabilistic sophistication and variational preferences.
\textit{Journal of Economic Theory} \textbf{146}, 2117--2125.


\bibitem[\protect\citeauthoryear{Tversky and Kahneman}{1992}]{TverskyKahneman92}
\textsc{Tversky, A., and Kahneman, D.} (1992).
Advances in prospect theory: Cumulative representation of uncertainty.
\textit{Journal of Risk and Uncertainty} \textbf{5}, 297--323.

\bibitem[\protect\citeauthoryear{Von Neumann and Morgenstern}{1944}]{NeumannMorgenstern44}
\textsc{Von Neumann, J., and Morgenstern, O.} (1944).
\textit{Theory of Games and Economic Behavior}.
3rd ed. 1953, Princeton University Press, Princeton.

\bibitem[\protect\citeauthoryear{Wakker}{2010}]{Wakker10}
\textsc{Wakker, P. P.} (2010).
\textit{Prospect Theory for Risk and Ambiguity}.
Cambridge University Press, Cambridge.

\bibitem[\protect\citeauthoryear{Wald}{1950}]{Wald50}
\textsc{Wald, A.} (1950).
\textit{Statistical Decision Functions}.
Wiley, New York.




\bibitem[\protect\citeauthoryear{Wang}{2022}]{Wang22}
\textsc{Wang, F.} (2022). Rank-dependent utility under multiple priors. 
\textit{Management Science} \textbf{68}, 8166--8183.

\bibitem[\protect\citeauthoryear{Wozabal}{2014}]{W14}
\textsc{Wozabal, D.} (2014).
Robustifying convex risk measures for linear portfolios: A nonparametric approach.
\textit{Operations Research} \textbf{62}, 1302--1315.

\bibitem[\protect\citeauthoryear{Yaari}{1969}]{Yaari69}
\textsc{Yaari, M. E.} (1969).
Some remarks on measures of risk aversion and on their uses.
\textit{Journal of Economic Theory} \textbf{1}, 315--329.


\bibitem[\protect\citeauthoryear{Yaari}{1987}]{Yaari87}
\textsc{Yaari, M. E.} (1987).
The dual theory of choice under risk.
\textit{Econometrica} \textbf{55}, 95--115.

\end{thebibliography}}



\end{document}